\documentclass[12pt]{article}
\usepackage{mathrsfs}
\usepackage{amsmath,amsfonts,amssymb,rotating,amsthm, bm}
\usepackage{hyperref}
\usepackage{array}
\usepackage{breqn}
\usepackage{color}
\usepackage{fancyheadings}
\usepackage{mathptmx}
\usepackage{enumitem}
\usepackage{dsfont}
\usepackage{caption}
\usepackage{subcaption}
\usepackage{booktabs}
\usepackage{tikz}
\usetikzlibrary{calc}
\usetikzlibrary{trees}

\tikzstyle{every node}=[circle,inner sep=1pt,fill=white!60]

\tikzstyle{tn}=[shape=circle, draw, color=black!70]
\tikzstyle{marke}=[shape=circle,minimum size=0.2cm, draw,blue]

\usepackage{times}

\allowdisplaybreaks
\def\qed{\nopagebreak\hfill{\rule{4pt}{7pt}}}

\newtheorem{thm}{Theorem}[section]

\newtheorem{cor}[thm]{Corollary}

\theoremstyle{remark}

\numberwithin{equation}{section}

\usetikzlibrary{trees}

\tikzstyle{every node}=[circle,inner sep=1pt,fill=white!60]
\tikzstyle{tn}=[shape=circle, draw, color=black!70]
\tikzstyle{tn1}=[shape=circle, draw]
\tikzstyle{tn2}=[shape=circle, draw,inner sep=1.5pt]
\tikzstyle{tn3}=[shape=rectangle, draw,inner sep=1.5pt]
\tikzstyle{marke}=[shape=circle,minimum size=0.1cm, draw,blue]

\newcommand\gen{{\rm Gen}}

\begin{document}

\begin{center}
{\large\bf A Grammatical Calculus for Peaks and Runs of Permutations}

\vskip 6mm

William Y.C. Chen$^1$ and Amy M. Fu$^2$

\vskip 3mm

$^{1}$Center for Applied Mathematics\\
Tianjin University\\
Tianjin 300072, P.R. China

\vskip 3mm

$^{2}$School of Mathematics\\
Shanghai University of Finance and Economics\\
Shanghai 200433, P.R. China

\vskip 3mm

Emails: { $^1$chenyc@tju.edu.cn, $^{2}$fu.mei@mail.shufe.edu.cn}

\vskip 6mm

    {\bf Abstract}
\end{center}

 We develop a nonstandard
approach  to exploring polynomials
associated with peaks and runs of permutations.
With the aid of a context-free grammar,
or a set of substitution rules,
one can perform a symbolic calculus, and the
computation  often  becomes rather simple.
From a grammar it follows at once a system of
ordinary differential equations for the generating functions.
Utilizing a certain constant property, it is even
possible to deduce a single equation for each generating function.
To bring the grammar to
a combinatorial setting,
we find a  labeling scheme for
up-down runs of a permutation,
which  can be regarded as a refined
 property, or the differentiability in a certain sense,
in contrast to the
usual counting argument for the recurrence relation. The labeling
scheme also exhibits how the substitution rules arise in the
construction of the combinatorial structures.
Consequently,
 polynomials on peaks and runs
can be dealt with in two ways,
combinatorially or  grammatically.
The grammar also serves as a guideline to
build a bijection between permutations and
increasing trees that maps the number of
up-down runs to the number of nonroot vertices of even degree.
This correspondence can be adapted to left peaks and exterior peaks, and the key
step of the construction is called the reflection principle.



\section{Introduction}

We present a grammatical approach to
studying peaks and runs of permutations.
This is a classical topic in enumerative
combinatorics, which has been extensively
investigated for decades, see B\'ona \cite{Bona-2004}. By injecting labels
or variables into the conventional procedure
of deriving recurrence relations, one obtains
a set of substitution rules, or a context-free
grammar. A grammar is more than just a recurrence relation,
it offers a ground for a grammatical calculus,
or a symbolic calculus in a rigorous sense, which
can be quite efficient in acquiring generating functions.
Besides, a grammar may shed light on the  construction of
a correspondence between two objects that have a grammar in
common. Armed with the grammar, we find
a bijection between permutations
and increasing trees that maps the number of
up-down runs and the number of peaks to the number of vertices
subject to certain degree conditions.

To be more specific, for $n\geq 1$, let $\sigma=\sigma_1\sigma_2\cdots \sigma_n$
be a permutation of $[n]=\{1,2,\ldots, n\}$. A left peak of $\sigma$ is
an index $i$ such that $1\leq i \leq n-1$ and $\sigma_{i-1} < \sigma_i >
\sigma_{i+1}$,
under the convention that $\sigma_0=0$.
An interior peak of $\sigma$ is an index $i$ such that
$2 \leq i \leq n-1$ and $\sigma_{i-1} < \sigma_i > \sigma_{i+1}$.
An exterior peak or an outer peak of $\sigma$ is an index
$1\leq i \leq n$ such that $\sigma_{i-1} < \sigma_i > \sigma_{i+1}$
with the understanding that $\sigma_0 = \sigma_{n+1}=0$.

For $n\geq 0$, the left peak polynomials, interior peak polynomials and
the exterior peak polynomials are denoted by $L_n(x)$, $M_n(x)$ and $W_n(x)$,
respectively.  Let $L(x,t)$,  $M(x,t)$ and $W(x,t)$ be the
exponential generating functions of $L_n(x)$, $M_n(x)$ and $W_n(x)$,
respectively.
 David-Barton \cite{David-Barton-1962}
 established partial differential equations
 for $L(x,t)$ and $M(x,t)$, and found a solution  requiring
 one more step of integration.
 An explicit expression of $L(x,t)$
was given by Gessel, see   \cite{Chen-Fu-2017, Chen-Fu-2022T, Zhuang-2016}.
A formula for $M(x,t)$ can also be deduced from
a generating function of  Carlitz-Scoville \cite{Carlitz-Scoville-1972}. Note
that $W(x,t)$ can be easily deduced from $M(x,t)$.
Alternative proofs of the formulas for $L(x,t)$ and $M(x,t)$
have been given by Zhuang \cite{Zhuang-2016}.

The grammar for the left peak polynomials $L_n(x)$ was discovered
by Ma \cite{Ma-2012} via a recurrence relation and  independently by Chen-Fu \cite{Chen-Fu-2017}
with a grammatical labeling.
Ma \cite{Ma-2012} further noticed that this grammar can be employed to generate
the interior peak polynomials.
Below is the concerned grammar
\begin{equation} \label{Gxy}
    G=\{ x \rightarrow xy, \;\, y\rightarrow x^2\}.
\end{equation}
Let $D$ denote the formal derivative with respect to $G$, or equivalently,
$D$ can be perceived as a differential operator
\begin{equation}
    D=xy {\partial \over \partial x} + x^2 {\partial \over \partial y}.
\end{equation}
Notice that the operator $D$ is a derivative so that the Leibniz formula
holds, and this property makes it possible to perform the grammatical
calculus with combinatorial motivations, see \cite{Chen-1993, Chen-Fu-2017, Dumont-1996}.
As remarked by Dumont \cite{Dumont-1996}, there are advantages to view $G$
as a set of substitution rules. In fact, the setting of a grammar
captures the combinatorial significance of a recursive procedure.

Despite previous attempts to explore the peak
polynomials by means of a grammar,
the story does not seem to be complete. For example, the authors
were nearly there in obtaining a direct derivation of the generating functions
solely by  conducting the grammatical computation, but we missed the ultimate
goal of reaching a sound understanding of the hyperbolic functions
appearing in the formulas. This regret or oversight is now redeemed by
considering the inverse of a generating function and by exploiting
a constant property in a grammatical sense.
It turns out that the constant property often facilitates the computation of a generating function.

One feature of the grammatical calculus is that
we need to work with
the bivariate versions of the peak polynomials. As far as the
coefficients are concerned, the variable $y$ is only a figurehead.
But it is just the opposite that in the grammatical calculus
the variable $y$ in its
own part is  as equally important as $x$. In fact, with the
company of the
 grammar, it is a pleasant journey to reach  bivariate versions of the
generating functions for the three kinds of peak polynomials.
Meanwhile, this explains why the hyperbolic functions crop up with no need to mention
differential equations.

As a bonus of the grammatical calculus, we get a simple
relationship between the generating functions of the left peak
polynomials and the exterior peak polynomials. In the context of peak
polynomials, we illustrate how a grammar can be automatically
translated into a system of
ordinary differential equations of the generating functions.
In light of a certain constant property, we may
even deduce a single ordinary differential equation for each generating function.

Then we move on to the number of up-down runs of
a permutation. The peaks and runs are closely related objects.
For a permutation $\sigma$ of $[n]$, we assume again that
a zero is patched at the beginning, i.e., $\sigma_0=0$.
An up-down run of $\sigma$ is
a maximal segment (a subsequence consisting  of consecutive elements)
that is either increasing or decreasing. If we drop the assumption
$\sigma_0=0$, then a maximal increasing or decreasing segment is called
an alternating run.
For example, the permutation $3\,  7\,  5\,  8\,  6\, 1\,   4\,  9\,  2 $
has six up-down runs:
\[ 0\  3\  7\ , \; 7\  5, \; 5\  8, \; 8\  6\ 1, \; 1\   4\  9, \; 9\  2\]
and it also has six alternating runs. However, the permutation
$2\,1$ has two up-down runs and only one alternating run.

For $n\geq 0$, let $\Lambda_n(x)$ denote the polynomial
associated with the number  of permutations of $[n]$
with $k$ up-down runs, and let $\Lambda(x,t)$ be
exponential generating function of $\Lambda_n(x)$.
Here we set our eyes on the  Greek letter $\Lambda$  because
it bears an ideal resemblance to the up-down shape and
we do not have to worry about interfering with the customary notation
$A_n(x)$ for the Eulerian polynomials.
As observed by B\'ona, $\Lambda(x,t)$ can be deduced from
a formula of David-Barton on the polynomials associated
with the number of permutations of $[n]$ with $k$ alternating runs.
An equivalent formula for $\Lambda(x,t)$  was obtained by
Stanley \cite{Stanley-2008} in regard with
permutations of $[n]$ having a given length of
the longest alternating subsequences.

The grammar that governs the counting of up-down runs was found by
Ma \cite{Ma-2013}, that is,
\begin{equation} \label{G-axy}
    G=\{a \rightarrow ax, \;\, x \rightarrow xy, \;\, y\rightarrow x^2\},
\end{equation}
which is equivalent to the operator
\begin{equation} \label{D-partial}
D = ax {\partial \over \partial a} +
     xy {\partial \over \partial x} +
      x^2 {\partial \over \partial y}.
\end{equation}
In the framework of the grammatical calculus,
we are led directly to a relationship between
$\Lambda_n(x)$ and the peak polynomials. We notice that
this connection also follows from two
relations due to Ma \cite{Ma-2013}.
Moreover, we exhibit that
the grammar is informative for deriving
exponential formulas on peaks and up-down
runs of permutations.

If one asks why there exists a
context-free grammar for the enumeration of
a combinatorial structure, a vague answer would be
that the involved quantities amount to local properties
with respect to the growth of the structure. Indeed,
this is what a context-free grammar is all about.
For example, when generating permutations of $[n+1]$
out of permutations of $[n]$, the operation of inserting
$n+1$ has rather a local impact on the number of various
peaks. Thus in a certain sense, a grammar is a way of
showing locality or differentiability.

For up-down runs, we give an explicit labeling of permutations,
called the up-down labeling, which reveals how the
insertion operation affects the number of up-down runs during
the procedure of generating a permutation of $[n+1]$ from
a permutation of $[n]$. With this labeling scheme in hand,
we see that the grammar is more than just a recurrence relation,
and we can make use of it in two ways, either as
 an apparatus of a grammatical calculus, or as a bridge
to facilitate finding
bijections between two objects that have a grammar in common.

This line of thinking yields a bijection between permutations
and increasing  trees that maps the number of up-down runs of a
permutation to the number of nonroot vertices of even degree
of an increasing tree. The grammar plays a vital role in
the justification of the bijection. For this reason, we call
the bijection a grammar assisted bijection.
In particular, when restricted to
alternating permutations, we arrive at a grammar assisted
correspondence
between alternating permutations and increasing even
trees. A construction in this regard has been given
by Kuznetsov, Pak and Postnikov \cite{KPP-1994}.

As for exterior  peaks, it is known that
the classical bijection between permutations and
increasing binary trees maps the number of the exterior peaks
to the number of certain type of vertices, see Stanley \cite{Stanley-V-1}.
Stemming from two
labeling schemes relative to the same grammar, we obtain
a bijection
between permutations and
increasing trees connecting the number of exterior
peaks of a permutation to the number of vertices
of even degree of an increasing tree.

It is our hope that this grammatical proposal could be   applicable
to more occasions.

\section{A grammatical calculus for peaks}

The objective of this section is to demonstrate the efficiency of the
grammatical calculus in the study of the three kinds of
peak polynomials. All we need
is a grammar and the Leibniz rule relative to the grammar.

 For $n\geq 1$ and $0\leq k \leq \lfloor n /2\rfloor$, let
 $L(n,k)$ denote the number of permutations of $[n]$ with
 $k$ left peaks. Analogously, $M(n,k)$ and $W(n,k)$ are
 defined in regard with interior peaks and exterior peaks,
 respectively, with $k$ in the valid range.
For $n=0$, we set $L(0,0)=M(0,0)=W(0,0)=1$. Define
\begin{eqnarray}
    L_n(x) & = & \sum_{k=0}^{\lfloor n /2\rfloor}
                  L(n,k) x^k, \\[6pt]
    M_n(x) & = & \sum_{k=0}^{\lfloor (n-1) / 2 \rfloor}
           M(n,k)x^k,\\[6pt]
      W_n(x) & = & \sum_{k=1}^{\lfloor (n+1) / 2 \rfloor}
           W(n,k)x^k.
\end{eqnarray}

In connection with the
grammar
\begin{equation} \label{Gx2y}
    G = \{  x \rightarrow xy, \;\;
           y \rightarrow x^2 \},
\end{equation}
for $n\geq 1$, the bivariate peak polynomials are defined by
\begin{eqnarray} \label{LMW-xy}
L_n(x,y) & = & \sum_{k=0}^{\lfloor n /2 \rfloor} L(n,k) x^{2k+1}
        y^{n-2k}, \\[6pt]
 M_n(x,y) & = &  \sum_{k=0} ^{\lfloor (n-1)/2\rfloor} M(n,k) x^{2k+2} y ^{n-2k-1}, \\[6pt]
 W_n(x,y) & = &  \sum_{k=1} ^{\lfloor (n+1)/2\rfloor} W(n,k)x^{2k} y ^{n-2k+1}.
\end{eqnarray}
For $n=0$, we define $L_0(x,y)= x$, $M_0(x,y)=1$ and $W_0(x,y)=y$. Let $D$ be the formal derivative with respect to $G$ given in \eqref{Gx2y}. We have the grammatical interpretations of the three bivariate peak polynomials, for $n \geq 0$,
$$
D^n(x)=L_n(x,y), \quad D^n(y)=W_n(x,y)
$$
and, for $n \geq 1$, $D^n(y)=M_n(x,y)$, respectively.

To perform the grammatical calculus, we need to consider the
generating function of  a Laurent polynomial
$f$ in $x$ and $y$  with respect to the operator $D$, as defined by
\begin{equation}
    \gen(f, t) = \sum_{n=0}^\infty D^n(f) {t^n \over n!}.
\end{equation}
 The formal derivative $D$ is related to the derivative
with respect to the variable $t$ of the generating function $\gen(f,t)$
via the following relation:
\begin{equation} \label{br}
\gen(D(f), t) = \gen'(f,t),
\end{equation}
where the prime signifies the derivative relative to $t$.

To see how the grammatical approach works, we  give a simple relation for the generating functions $L(x,y,t)$ and $W(x,y,t)$, that is,
$$
L(x,y,t)=  \gen(x, t),  \quad  W(x,y,t)=  \gen(y, t).
$$
Then we present a combinatorial interpretation of this fact.

\begin{thm} \label{thm-LW}
We have
\begin{equation}\label{LW-1}
    \gen^2(x,t) = \gen^2(y,t) + x^2 -y^2.
\end{equation}
\end{thm}

\noindent
{\it Proof.} It is easily checked that
$x^2-y^2$ is a constant relative to $D$,
that is, $D(x^2) = D(y^2)$, see \cite{Chen-Fu-2017}.
Thus we have
$D^n(x^2)=D^n(y^2)$ for $n\geq 1$, and hence
\begin{equation}
    \gen(x^2,t) =x^2+ \sum_{n=1}^\infty D^n(x^2){t^n \over n!}=
   x^2+ \sum_{n=1}^\infty D^n(y^2){t^n \over n!},
\end{equation}
which is the right-hand side of (\ref{LW-1}), as claimed.  \qed

Unlike the approach of David-Barton by means
of partial differential equations,
in virtue of the above relation (\ref{LW-1}), we may
aim at  ordinary
differential equations for $L(x,y,t)$ and $W(x,y,t)$, which are within the reach of
Maple. Note that we need the constant property in deriving the equations. This
way of deriving ordinary differential equations may suit other instances.
 For example, for the grammar $G=\{x\rightarrow xy, \;\, y\rightarrow xy\}$, which generates the Eulerian polynomials, $x-y$ is a constant. For the grammar
$ G=\{x \rightarrow xy, \;\, y \rightarrow x\}$, which generates
the Andr\'e polynomials, $y^2-2x$ is a constant.  In each case,
one can deduce
 an ordinary differential equation from the grammar.
 In fact, once an equation for $\gen(x,t)$
 is obtained,
 we may treat $y$ merely as a parameter and set $y=1$ for the purpose of
 computing the generating function only involving $x$.
 We also remark that the relation (\ref{br}) enables
 us to read off a system
 of ordinary differential equations on the generating functions
 for all the variables. As for the peak polynomials,  we come to
 the following system of equations reminiscent of the grammar,
 \begin{equation}
         \left\{  \begin{aligned}   L'( t) & =   L (t) W(t), \\[3pt]
              W'(t)  & =   L^2(t),\\[1pt]
                  \end{aligned}
              \right.
 \end{equation}
 with boundary conditions $L(0)=x$ and $W(0) = y$,
where the parameters $x$ and $y$ in $L(x,y,t)$ and $W(x,y,t)$ are suppressed to emphasize that the
derivative is taken with respect to $t$.

\begin{thm} The following ordinary differential equations hold
with
boundary conditions $L(x,y,0)=x$ and $W(x,y,0)=y$:
\begin{eqnarray}
    L'(t) & = & L(t) \sqrt{ L^2(t) - x^2+y^2},\\[3pt]
    W'(t) & = & W^2(t)+x^2-y^2.
\end{eqnarray}
\end{thm}

There is a combinatorial explanation of the fact that
$D^{n}(x^2) = D^n (y^2)$ for $n\geq 1$. Since $D(y)=x^2$,
$D^n(x^2)$ can be rewritten as $D^{n+1}(y)$. Thus, the combinatorial reason behind the relation (\ref{LW-1})
lies in the following two convolution formulas for $n\geq 1$,
\begin{eqnarray}
W_{n+1}(x,y) &  = & \sum_{k=0}^n {n \choose k} W_{k}(x,y) W_{n-k}(x,y), \label{WW}\\[6pt]
W_{n+1}(x,y) &  = & \sum_{k=0}^n {n \choose k} L_k(x,y) L_{n-k}(x,y).  \label{WL}
\end{eqnarray}

For a permutation $\sigma$, we use
$L(\sigma)$, $M(\sigma)$ and $W(\sigma)$ to
denote the number of left peaks, interior peaks and
exterior peaks of $\sigma$, respectively.
A bijective argument for (\ref{WW}) and (\ref{WL}) goes as follows.
For $n\geq 1$, let $\sigma$ be a permutation of $[n+1]$.
Write $\sigma= \pi 1 \tau$. It is readily seen
\[ W(\sigma) = W(\pi) + W(\tau).\]
Thus we arrive at (\ref{WW}).

On the other hand, let us write
$\sigma=\pi (n+1) \tau$, and let $\tau'$ denote the
reverse of $\tau$.
We find that
\[ W(\sigma) = L(\pi) + L(\tau') + 1,\]
since $n+1$ is always an exterior peak.
This special exterior peak is counted by the
exponent of $x$ in the definition of $L_n(x,y)$, with
each contributing a factor $x$.

With the above two interpretations of $W_{n+1}(x,y)$, we
may burn the bridge after having crossed the river.
In doing so, we get
a direct correspondence between the right-hand sides
of (\ref{WW}) and (\ref{WL}).

 The next theorem gives the generating function of $x^{-1}$
with respect to $D$, which is the inverse of the bivariate form of Gessel's formula.

\begin{thm} We have
\begin{equation} \label{bg}
    \gen(x^{-1},t)  =\frac{\sqrt{y^2-x^2}\cosh(t\sqrt{y^2-x^2})-y\sinh(t\sqrt{y^2-x^2})}
    {x\sqrt{y^2-x^2}}.
\end{equation}
\end{thm}

\noindent
{\it Proof.}
As noted in \cite{Chen-Fu-2017},
\[
D(x^{-1}) = -x^{-1}y, \quad
D^2(x^{-1}) = x^{-1}(y^2-x^2).
\]
Now the following pattern emerges. For $n\geq 0$,
\begin{eqnarray}
D^{2n}(x^{-1})& =& x^{-1}  (y^2-x^2)^n, \label{x2n}\\[6pt]
 D^{2n+1}(x^{-1}) &= &-x^{-1}y(y^2-x^2)^n. \label{x2n1}
\end{eqnarray}
Taking the parity into account,  we find that
\begin{eqnarray}
\sum_{n=0}^\infty D^{2n}(x^{-1}) {t^{2n} \over (2n)!} & =& x^{-1}\cosh(t\sqrt{y^2-x^2}),\\[6pt]
 \sum_{n=0}^\infty D^{2n+1}(x^{-1})
{t^{2n+1} \over  (2n+1)!} & = & -\frac{x^{-1}y}{\sqrt{y^2-x^2}}\sinh(t\sqrt{y^2-x^2}).
\end{eqnarray}
Putting the above sums together completes the proof.
\qed

\begin{cor} We have
\begin{equation}\label{Genx}
    \gen(x,t)  =\frac{x\sqrt{y^2-x^2}}
    {\sqrt{y^2-x^2}\cosh(t\sqrt{y^2-x^2})-y\sinh(t\sqrt{y^2-x^2})}.
\end{equation}
\end{cor}

Let $L(x,t)$ be the exponential generating function of $L_n(x)$. Dividing both sides
by $x$, replacing $x^2$ by $x$ and setting $y=1$ yields
Gessel's formula:
\begin{equation}
    L (x,t)  = { \sqrt{1-x} \over
      \sqrt{1-x} \cosh  (t \sqrt{1-x})
              - \sinh (t \sqrt{1-x}) }.
\end{equation}

Given the above formula for $\gen(x, t)$,
the  generating function $\gen(y,t)$ for exterior  peaks
can be deduced from the relation
$$
\gen'(x,t)=\gen(x,t)\gen(y,t),
$$
or equivalently,
\begin{equation}\label{xy}
\gen(y,t)=\frac{\gen'(x,t)}{\gen(x,t)}.
\end{equation}

\begin{cor} We have
\begin{equation} \label{Geny}
\gen(y,t) =
\frac{y\sqrt{y^2-x^2}\cosh(t\sqrt{y^2-x^2})-(y^2-x^2)\sinh(t\sqrt{y^2-x^2})}
{\sqrt{y^2-x^2}\cosh(t\sqrt{y^2-x^2})-y\sinh(t\sqrt{y^2-x^2})}.
\end{equation}
\end{cor}

There are alternative ways to derive the above formula for
$\gen(y,t)$. Analogous to $\gen(x^{-1},t)$, it is easy to compute
$\gen(x^{-1}y, t)$.

\begin{thm} We have
\begin{equation}\label{x-1-y-3}
    \gen (x^{-1}y, t)   =    x^{-1} y\cosh(t \sqrt{y^2-x^2})
       -x^{-1}\sqrt{y^2-x^2}\sinh(t\sqrt{y^2-x^2}).
\end{equation}
\end{thm}

Notice that the above relation also implies
the formula (\ref{Geny}) for $\gen(y,t)$, since
\begin{equation} \label{yxt}
\gen(y,t) = \gen(x,t) \gen(x^{-1}y, t) .
\end{equation}
One more  way to relate $\gen(y,t)$ to
$\gen(x,t)$ is to utilize the relation (\ref{LW-1}).

Let $M(x,y,t)$ be the exponential generating function of $M_n(x,y)$. Owing to the fact
\[
D^n(y)= M_n(x,y) ,
\]
for $n\geq 1$,  and the initial value $M_0(x,y)=1$,
we see that
\begin{equation}
\label{M-y}
M(x,y,t)=1-y+ \gen(y,t),
\end{equation}
 which implies the David-Barton formula in the bivariate form.

 \begin{cor} We have
 \begin{equation}
\sum_{n=0}^\infty M_n(x,y)\frac{t^n}{n!}=\frac{\sqrt{y^2-x^2}\cosh(t\sqrt{y^2-x^2})+(1-y)\sinh(t\sqrt{y^2-x^2})}{\sqrt{y^2-x^2}\cosh(t\sqrt{y^2-x^2})-y\sinh(t\sqrt{y^2-x^2})}.
\end{equation}
\end{cor}

Let $M(x,t)$ be the exponential generating function of $M_n(x)$. Setting $y=1$  and replacing $x^2$ by $x$, we get
\begin{equation}
    M(x,t)  =
    { \sqrt{1-x} \cosh (t \sqrt{1-x})
      \over
    \sqrt{1-x} \cosh (t \sqrt{1-x})
              - \sinh (t \sqrt{1-x})} .
\end{equation}

Considering a grammar
\begin{equation}
    H=\{ a \rightarrow ay, \;\, x\rightarrow xy, \;\, y \rightarrow x^2\},
\end{equation}  we find that $L_n(x)$ and $W_n(x)$ satisfy an
exponential relation.
It is not hard to see that
$D^n(a)$ relative to $H$ coincides with
$D^n(x)$ relative to  $G=\{x\rightarrow xy, \;\,y\rightarrow x^2\}$.

\begin{thm} \label{thm-LW-E} We have
\begin{equation} \label{GLW-3}
  \sum_{n=0}^\infty L_{n}(x) {t^n \over n!}
  = \exp \left( \sum_{n= 0}^\infty  W_{n }(x) {t^{n+1} \over (n+1)!} \right).
\end{equation}
\end{thm}

Next, we seek a combinatorial interpretation
of  (\ref{GLW-3}). To this end, we
introduce a decomposition
of a permutation, called the $LW$-decomposition, which is
essentially the cycle decomposition, or Foata's first fundamental transformation of a permutation.
Assume that $\sigma$ is
a permutation of $[n]$. If $1$ appears at the end of $\sigma$, then
nothing needs to be done. Otherwise, write $\sigma=\pi  \tau$, where
$\pi$ ends with $1$. Now we make $\pi$ the first block  of the
decomposition and repeat the process for $\tau$. If $\tau$
ends with the minimum element, then we are done. Otherwise, we continue to
decompose $\tau$ in the same manner.
For example, the permutation $2\, 6\, 1\, 3\,  8\, 4\,
7\, 9\,  5$ is decomposed into four segments:
\[  2\  6\ 1 \  | \  3 \ | \  8\   4 \  | \  7\ 9\  5  .\]
Notice that the blocks are displayed in the increasing order of their  minimum
elements.

\noindent
{\it Proof of Theorem \ref{thm-LW-E}.}
Let $\sigma$ be a permutation of $[n]$, where $n\geq 1$, and
let $\sigma=\pi_1|\pi_2|\cdots |\pi_k$ be the $LW$-decomposition
of $\sigma$. For $1\leq i \leq k$, let $\pi'_i$
be the permutation obtained from $\pi_i$ by removing
the minimum element at the end. We proceed to show that
\begin{equation} \label{LW-sp}
    L(\sigma) = W(\pi'_1 ) + W(\pi'_2 ) + \cdots + W(\pi'_k ).
\end{equation}
Asssume that
the minimum element of $\pi_i$ is the element $\sigma_j$.
We claim that $j$ cannot be a left peak in $\sigma$.
Otherwise, the element $\sigma_{j+1}$ would be smaller than $\sigma_{j}$, which is
contradictory to the construction of the $W$-decomposition.

On the other hand, assume that $\sigma_j$ is not at
the end of any segment $\pi_i$ of the $LW$-decomposition
of $\sigma$, say $\sigma_j$ is an element of $\pi_i'$.
It is evident that
$j$ is a left peak of $\sigma$ if and only if the corresponding position
is an exterior peak of $\pi_i'$. Thus we arrive at (\ref{LW-sp}).
This completes the proof.
\qed

\section{A grammatical calculus for up-down runs}

Up-down runs of a permutation are intuitively related to
peaks. It is also known that
the number of up-down runs of a permutation of $[n]$ equals
the length of the longest alternating
subsequences, see Stanley \cite{Stanley-2008}. From the viewpoint
of grammars, the structure of up-down runs can be understood as
an exponential structure built on left peaks. This
connection will be made precise in Theorem \ref{thm-AL}.
This section is devoted to a grammatical calculus for the
number of permutations of $[n]$ with $k$ up-down runs.

For $n\geq 2$, let $\Lambda(n,k)$ denote the number of permutations   of $[n]$ with
$k$ up-down runs. In particular, define $\Lambda(0,0)=\Lambda(1,1)=1$ and
$\Lambda(0,k)=\Lambda(n,0)=0$ whenever $n,k\geq 1$.  For $n\geq 0$,
write
\begin{equation}
    \Lambda_n(x) =\sum_{k=1}^n \Lambda(n,k) x^k.
\end{equation}

The following grammar was discovered by Ma \cite{Ma-2013}:
\begin{equation}
    G=\{a \rightarrow ax, \;\, x\rightarrow xy, \;\, y\rightarrow x^2\}.
\end{equation}
For $n\geq 0$, the bivariate form of $\Lambda_n(x)$
is defined by
\begin{equation}\label{A-nxy}
    \Lambda_n(x,y)= \sum_{k=1}^{n} \Lambda(n,k) x^k y^{n-k}.
\end{equation}
In view of the recurrence relation
\begin{equation} \label{an-r}
  \Lambda(n,k)= k \Lambda(n-1,k) + \Lambda(n-1,k-1) + (n-k+1) \Lambda(n-1,k-2),
\end{equation}
  for $n,k\geq 1$, Ma obtained the following expression
  for $\Lambda_n(x,y)$ in terms of the above grammar $G$.

  \begin{thm}[Ma \cite{Ma-2013}] \label{thm-ma}
  For $n\geq 0$, we have
\begin{equation} \label{Ma-Da}
    D^n(a) = a \Lambda_n(x,y).
\end{equation}
  \end{thm}

The first few values of $D^n(a)$ are given below:
\begin{eqnarray*}
D(a) & = & ax, \\[6pt]
D^2(a)  & = &  a xy+a x^2, \\[6pt]
D^3(a) & = &  a  xy^2+3 a  x^2y+2 a x^3 , \\[6pt]
D^4(a) & = & a  xy^3+7 a  x^2y^2+11 a x^3 y+5 a x^4,\\[6pt]
D^5(a) & = &  ax  y^4 +15ax^2 y^3+43 a x^3 y^2 +45 a x^4y +16 a x^5, \\[6pt]
D^6(a) & = & a  xy^5+31 a x^2 y^4 +148 a x^3y^3 +268 a x^4y^2 +211 a x^5y +61 a x^6 .
\end{eqnarray*}

As observed by B\'ona \cite{Bona-2004, Bona-2021}, see also  Stanley \cite{Stanley-2010},
$\Lambda_n(x)$ can be expressed in terms
of the polynomials $R_n(x)$ associated with alternating runs,
and in turn $R_n(x)$ can be expressed in terms of the Eulerian polynomials
as shown by David-Barton \cite{David-Barton-1962}, see also Knuth
\cite{Knuth}.
More precisely, for $n\geq 1$ and $1\leq k \leq n$,
let $R(n,k)$ denote the number of permutations of $[n]$ with $k$ alternating
runs.
For $n\geq 0$,
let
\begin{equation}
    R_n(x,y)=\sum_{k=1}^n R(n+1,k) x^k y^{n-k}.
\end{equation}
Ma \cite{Ma-2013} showed that for $n\geq 0$,
\begin{equation}
    D^n(a^2) = a^2 R_n(x,y).
\end{equation}
Let $R(x,t)$ denote the exponential generating function of $R_n(x,y)$ with
$y$ set to $1$, and let $\Lambda(x,t)$ denote the exponential generating function of $\Lambda_n(x)$. With respect to the grammar $G$,
the relation
\[ \gen(a^2,t) = \gen^2(a,t)\]
takes the form of the following identity on generating functions,
which seems to have been unnoticed before, at least not explicitly,
\begin{equation}\label{RA}
R(x,t)= \Lambda^2(x,t).
\end{equation}

The grammatical labeling to be given in the next section
can be regarded as a combinatorial interpretation
of Theorem \ref{thm-ma}.

To explore the connections between
$\Lambda_n(x)$ and the peak polynomials from the aspect of the
grammar $G$,  we present an exponential
formula for $\Lambda_n(x)$ and $L_n(x)$. Then we provide
a grammatical derivation of a relation on $\Lambda_n(x,y)$,
$L_n(x,y)$ and $W_n(x,y)$.
Moreover, we
find a transformation of grammars leading to a relation on
$\Lambda_n(x)$ and $M_n(x)$ due to Ma \cite{Ma-2013}.

\begin{thm} \label{thm-AL} We have
\begin{equation} \label{AL-1}
  \sum_{n=0}^\infty \Lambda_n(x) {t^n \over n!}
  = \exp \left( \sum_{n= 0}^\infty  x L_{n }(x^2) {t^{n+1} \over (n+1)!} \right).
\end{equation}
\end{thm}

To give a combinatorial interpretation of the above theorem, we
present a decomposition of a permutation $\sigma$ such that
the number of up-down runs of $\sigma$ can be computed from the
numbers of left peaks of its components. This decomposition is called the
$AL$-decomposition. Let $n\geq 0$ and let $\sigma$ be a permutation of $[n]$.
  First, we observe that the number of up-down runs of the permutation
$\pi=\sigma(n+1)$ equals $2L(\sigma)+1$.   Therefore, $xL_n(x^2)$ equals
the generating function of permutations of $[n+1]$ ending with the maximum
element $n+1$ with respect to the number of up-down runs.

Given a permutation $\sigma$ of $[n]$, we may write $\sigma=\pi_1 \pi_2 \cdots \pi_k$
such that $\pi_1$ ends with the maximum element of $\sigma$, $\pi_2$ ends with
the maximum element of $\pi_2\cdots \pi_k$, and so forth.
This decomposition can be considered as a combinatorial interpretation of (\ref{AL-1}).

The $LW$-decomposition and the $AL$-decomposition can be reckoned as
the dual of each other.
 The following
identity is a bivariate form of a relation
observed by Ma \cite{Ma-2013}.

 \begin{thm} \label{AL-xy} We have
 \begin{equation}
\label{gen-at}
{(x+y)}\gen(a,t) =a \left(\gen(x,t)+\gen(y,t)\right),
\end{equation}
 or equivalently, for $n\geq 0$,
\begin{equation} \label{A-L-W}
  {(x+y) }   \Lambda_n(x,y) = L_n(x,y) + W_n(x,y).
^{}\end{equation}
\end{thm}

As a matter of fact,  it is a direct
consequence of following relations also given by  Ma \cite{Ma-2013}. For $n\geq 1$ and
$1 \leq k \leq \lfloor (n-1)/2 \rfloor$, we have
\begin{equation} \label{LA}
L(n,k)  =  \Lambda(n,2k) + \Lambda(n,2k+1),
\end{equation}
and for $n\geq 1$ and $0 \leq k \leq \lfloor (n-2)/2\rfloor$,
\begin{equation} \label{MA}
M(n,k)   =   \Lambda(n,2k+1) + \Lambda(n,2k+2).
\end{equation}

We now proceed to  compute
$\gen(a, t)$ using the grammatical calculus.

\noindent
{\it Proof.}
Observe that
 $\frac{x+y}{a}$ is a constant with respect to grammar, that is,
$$
D\left(\frac{x+y}{a}\right)=0.
$$
Consequently,
$$
  \gen(x,t)+ \gen(y,t)=\gen(x+y, t)=\gen\left(a \cdot \frac{x+y}{a}, t\right)=\frac{x+y}{a}\gen(a,t),
$$
as claimed. \qed

Substituting \eqref{Genx} and \eqref{Geny} into (\ref{gen-at}) and using the
exponential forms for $\sinh$ and $\cosh$, we obtain the following formula.

\begin{cor} We have
\begin{equation}
\gen(a,t) = a(y-x)\frac{y+\sqrt{y^2-x^2}+2xe^{\sqrt{y^2-x^2}t}+(y-\sqrt{y^2-x^2})e^{2\sqrt{y^2-x^2}t}}{y^2-x^2+y\sqrt{y^2-x^2}+(y^2-x^2-y\sqrt{y^2-x^2})e^{2\sqrt{y^2-x^2}t}}.
\end{equation}
\end{cor}

Setting $a=y=1$ yields  Stanley's formula \cite{Stanley-2010}:
\begin{equation}\label{SF}
    \Lambda(x,t) = (1-x) { 1 + \rho +2xe^{\rho t} + (1-\rho ) e^{2\rho t}
    \over 1+ \rho -x^2 + (1-\rho -x^2) e^{2\rho t}},
\end{equation}
where $\rho = \sqrt{1-x^2}$.

We conclude this section with a grammatical derivation of the
following relation on $\Lambda_n(x)$ and $M_n(x)$ due to  Ma \cite{Ma-2012}.
 For $n\geq 0$,
\begin{equation}\label{AM}
    \Lambda_n(x) = { x (1+x)^{n-1} \over 2^{n-1}}
     M_n\left( {2x \over 1+x} \right).
    \end{equation}

\noindent
{\it Grammatical Proof of \eqref{AM}.}
Set
\[ u=x+y,\;\,   v=\sqrt{x(x+y)}.\]
The grammar
$G=\{
x \rightarrow xy, \,\; y \rightarrow x^2
\} $
is transformed into
$$ G=\left\{
u \rightarrow v^2, \,\; v \rightarrow \frac{uv}{2}\right\}.
$$
Using the grammatical interpretation of the interior peak
polynomials $M_n(x,y)$, we find that
\begin{equation}\label{UV}
D^n(u)=v^2 \sum_{k=0}^{\lfloor (n-1)/2 \rfloor }
{1 \over 2^{n-1-k}} \, M(n,k) u^{n-1-2k}v^{2k}.
\end{equation}
By Theorem \ref{AL-xy}, we see that
\begin{equation}\label{DAM}
D^n(a)=\frac{a}{x+y}D^n(x+y),
\end{equation}
which yields
\begin{eqnarray*}
D^n(a)&=&ax \sum_{k=0}^{\lfloor (n-1)/2 \rfloor }{1 \over 2^{n-1-k}} \, M(n,k)  (x+y)^{n-1-2k}(x(x+y))^k\\[6pt]
&=& ax \frac{(x+y)^{n-1}}{2^{n-1}}\sum_{k=0}^{\lfloor (n-1)/2 \rfloor } M(n,k)\frac{(2x)^k}{(x+y)^k}.
\end{eqnarray*}
Setting $a=y=1$ completes the proof. \qed

\section{A grammatical labeling for up-down runs}

For a permutation $\sigma$ of $[n]$, the up-down labeling
of $\sigma$ can be described as follows.
Like the $L$-labeling given in \cite{Chen-Fu-2022T}, the labels are assigned
to the positions next to each element of $\sigma$.
For $1\leq i \leq n+1$, by a position $i$ we mean the position
immediately before $\sigma_i$, whereas the position $n+1$ is meant
to be the position after $\sigma_n$. Associated with
 a labeling, the weight of a permutation is
referred to the product of the labels.

To define the labeling for up-down runs of a permutation $\sigma$,
which we call the $A$-labeling, we first patch a zero
at the beginning of $\sigma$, or equivalently, set $\sigma_0=0$.
Then the labels are given by the following procedure.

\noindent
Case 1. For each up run
$\sigma_i \sigma_{i+1} \cdots \sigma_{j}$ possibly with $i=0$,
where $i<j\leq n$, two possibilities arise.
If $j<n$, that is, $\sigma_j$ is not the last element of $\sigma$,
we label the positions $i+1, \ldots, j-1$ by
$y$ and label the position $j$ by $x$.

If $j=n$, that is, $\sigma_j$ is the last element of $\sigma$,
then label the position $n$ by $a$ and position $n+1$ by $x$, and
the other positions by $y$. This case looks a little peculiar, but
that is perhaps the way it is.

\noindent
Case 2.
For each down run
$\sigma_i \sigma_{i+1} \cdots \sigma_{j}$, where $i<j$,
we always label the position $i+1$ by $x$,
and label the other positions $i+2, \ldots, j$ by $y$.
Moreover, if $j=n$, we label the position $n+1$ by $a$.

 For example,
below is the $A$-labeling of a permutation ending with a down run:
\[   0 \  y \     3 \   x  \  7 \  x \
  5  \   x \  8 \   x  \  6 \   y \   1 \ y\
 \ 4 \ x \ 9 \ x \  2 \  a.\]
The weight of the above permutation equals $ax^6y^3$.
For an example of $\sigma$ ending with
an up run, let $n=9$ and $\sigma= 3\,  7\,  5\,  8\,  6\,  1\,  2\,  4\,  9$.
Below is the
$A$-labeling,
\[   0 \  y \     3 \   x  \  7 \  x \
  5  \   x \  8 \   x  \  6 \   y \   1 \ y\
 \ 2 \  y \ 4 \ a \  9 \  x.\]
In this case the weight of $\sigma$  equals   $ax^5y^4$.

A labeling of a permutation is called a grammatical
labeling with respect to $G$, we mean that
when generating the permutations of $[n+1]$ by
inserting the element $n+1$ into a permutaiton
$\sigma$ of $[n]$, the substitution rules are applied
to each label exactly once.
This property makes it possible to
compute the total weight (sum of the products of labels)
of permutations of $[n+1]$ from the total weight
of permutations of $[n]$ by taking the formal derivative
with respect to the grammar. Moreover, the labels can be
treated as ingredients of the
 grammatical calculus. The notion of a grammatical labeling
was introduced in \cite{Chen-Fu-2017}, which  says that
the combinatorial structure is differentiable in a certain sense.
The following theorem
justifies the grammatical labeling for the purpose of
updating the weights of permutations upon the
insertion operations.

\begin{thm}\label{thm-ud}
 Let $n\geq 1$ and let $\sigma$ be a permutation of $[n]$.
Assume that $1 \leq i \leq n+1$.
Let $\pi$ be the permutation obtained from $\sigma$ by
inserting $n+1$ in $\sigma$ at  position $i$. Then
the weight of $\pi$ can be derived from that of $\sigma$ by
applying the substitution rule to the label of the element
of $\sigma$ at
position $i$.
\end{thm}

Before presenting the proof, let us give an  example.
Let $n=5$ and
$\sigma=2\,  5\,  4\,  1\,  3$.
In the table below,
an underlined label  signifies where the element $6$ is
inserted which is also the label to which
the substitution rule is applied.
\begin{center}
\renewcommand\arraystretch{1.33}
\begin{tabular}{| l|l|l|}
\hline
$\sigma$ & $\pi$ & Substitution\\ \hline
$0\  \underline{y}
 \ 2 \  {x} \  5 \  x \ 4 \  y\  1 \  a \ 3 \  x  $ \qquad\quad &
  $0 \  x \  6 \  x \   2 \   x  \   5 \  x  \   4 \  y \  1 \  a \   3 \   x $
       \qquad\quad  & $ y \rightarrow x^2$ \\ \hline
$0\  {y} \ 2 \  \underline{x} \  5 \  x \ 4 \  y\  1 \  a \ 3 \  x  $ \qquad&
  $0 \  y \   2 \ x\  6 \  x  \   5 \  y  \   4 \  y \  1 \  a \   3 \   x $
       \qquad & $ x \rightarrow xy$ \\ \hline
 $0\  {y} \ 2 \  {x} \  5 \  \underline{x} \ 4 \  y\  1 \  a \ 3 \  x  $ \qquad&
  $0 \  y \   2 \ y\  5 \  x  \ 6 \  x \  4 \  y \  1 \  a \   3 \   x $
       \qquad & $ x \rightarrow xy$ \\ \hline
 $0\  {y} \ 2 \  {x} \  5 \   {x} \ 4 \  \underline{y}\  1 \  a \ 3 \  x  $ \qquad&
  $0 \  y \   2 \ x\  5 \  x  \  4 \  x \ 6 \  x \  1 \  a \   3 \   x $
       \qquad & $ y \rightarrow x^2$ \\ \hline
$0\  {y} \ 2 \  {x} \  5 \   {x} \ 4 \   {y}\  1 \  \underline{a} \ 3 \  x  $ \qquad&
  $0 \  y \   2 \ x\  5 \  x  \  4 \  y \ 1 \ x\   6 \  x \    3 \   a $
       \qquad & $ a \rightarrow ax$ \\ \hline
$0\  {y} \ 2 \  {x} \  5 \   {x} \ 4 \   {y}\  1 \  {a} \ 3 \
\underline{x}  $ \qquad&
  $0 \  y \   2 \ x\  5 \  x  \  4 \  y \ 1 \ y\   3 \   a \  6 \  x \  $
       \qquad & $x   \rightarrow xy$ \\ \hline
  \end{tabular}
\end{center}

\noindent
{\it Proof of Theorem \ref{thm-ud}.}
It is readily seen that the Theorem is valid for $n=1,2$.
Now assume that
$n\geq 2$ and $\sigma$ is a permutation of $[n]$.
Suppose that $\pi$ is a permutation of $[n+1]$
created from $\sigma$ by
inserting the element $n+1$ at the position before $\sigma_i$, where
$1\leq i \leq n$, or at the position $n+1$, that is, at the end
of $\sigma$.

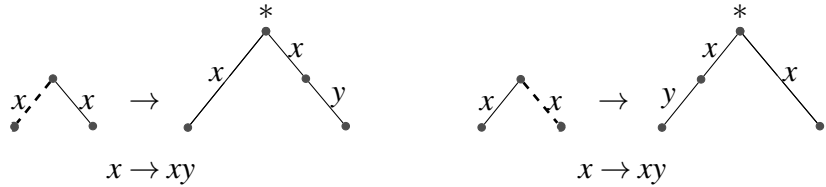
\begin{figure}[!ht]
\begin{minipage}{0.51\linewidth}
\begin{center}
\hspace*{1.3cm}
\begin{tikzpicture}[scale=0.55]
\node (zero)[tn,label=0:{}]{}[grow=down]
	%[sibling distance=16mm,level distance=10mm]
    child[grow=231,dashed,line width=1pt] {node [tn, label=0:{}](one){}}
    child[grow=310] {node [tn, label=0:{}](two){}
        };
    \coordinate[label=left:{$x$}] (a) at ($(zero)!.5!(one)$);
     \coordinate[label=right:{$x$}] (x) at ($(zero)!.5!(two)$);
\end{tikzpicture}
\hspace{1mm}
\raisebox{0.3cm}{$\rightarrow$}
\hspace{1mm}
\begin{tikzpicture}[scale=0.55]
\node (zero)[tn,label=0:{}]{}[grow=down]
	%[sibling distance=16mm,level distance=10mm]
    child[grow=231] {node [label=180:{}](one){}
     %[sibling distance=14mm,level distance=13mm]
            child {node [tn,label=0:{}](two){}
        }}
    child[grow=310] {node [tn,label=0:{}](three){}
    child[grow=310] {node [tn,label=0:{}](four){}}
        };
        \draw(zero)--(two);
       \node [tn,label=above:{$*$}] (1) at
(zero) {};
\coordinate[label=right:{$x$}] (x) at ($(zero)!.4!(three)$);
\coordinate[label=left:{$x$}] (Xx) at ($(zero)!.45!(two)$);
    \coordinate[label=right:{$y$}] (y) at ($(three)!.4!(four)$);
\end{tikzpicture}\\
\hspace*{0.6cm}$x\rightarrow xy$
\end{center}
\end{minipage}\hspace{-1.8cm}
\begin{minipage}{0.51\linewidth}
\begin{center}
\hspace*{1.4cm}
\begin{tikzpicture}[scale=0.55]
\node (zero)[tn,label=0:{}]{}[grow=down]
	%[sibling distance=16mm,level distance=10mm]
    child[grow=231] {node [tn, label=0:{}](one){}}
    child[grow=310,dashed,line width=1pt] {node [tn, label=0:{}](two){}
        };
    \coordinate[label=left:{$x$}] (a) at ($(zero)!.5!(one)$);
     \coordinate[label=right:{$x$}] (x) at ($(zero)!.5!(two)$);
\end{tikzpicture}
\hspace{1mm}
\raisebox{0.3cm}{$\rightarrow$}
\hspace{1mm}
\begin{tikzpicture}[scale=0.55]
\node (zero)[tn,label=0:{}]{}[grow=down]
	%[sibling distance=16mm,level distance=10mm]
    child[grow=231] {node [tn, label=180:{}](one){}
     %[sibling distance=14mm,level distance=13mm]
            child {node [tn,label=0:{}](two){}
        }}
    child[grow=310] {node [label=0:{}](three){}
    child[grow=310] {node [tn,label=0:{}](four){}}
        };
         \draw(zero)--(four);
       \node [tn,label=above:{$*$}] (1) at
(zero) {};
\coordinate[label=left:{$x$}] (x) at ($(zero)!.4!(one)$);
    \coordinate[label=left:{$y$}] (y) at ($(one)!.4!(two)$);
        \coordinate[label=right:{$x$}] (xy) at ($(zero)!.45!(four)$);
\end{tikzpicture}\\
\hspace*{0.7cm}$x\rightarrow xy$
\end{center}
\end{minipage}
\caption{Insertion at a position labeled by $x$.}
\label{fx}
\end{figure}

First, we consider the case when $i<n$.
If the position  is labeled by $x$,  it
can be seen that if $n+1$ is inserted in $\sigma$ at position $i$,
the change of weights is connected with the substitution
rule $x\rightarrow xy$. The two possibilities are illustrated in Figure \ref{fx},
where $*$ stands for the element $n+1$ and a dotted line indicates
the position of insertion.

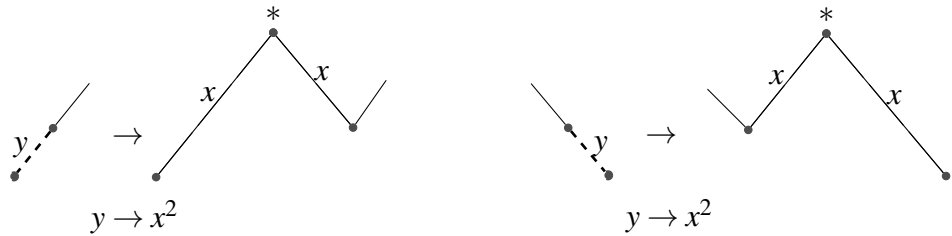
\begin{figure}[!ht]
\begin{minipage}{0.51\linewidth}
\begin{center}
\hspace*{1.3cm}
\begin{tikzpicture}[scale=0.55]
\node (zero)[label=0:{}]{}[grow=down]
	%[sibling distance=16mm,level distance=10mm]
    child[grow=231] {node [tn, label=0:{}](one){}
    child[dashed,line width=1pt] {node [tn, label=0:{}](two){}}
        };
    \coordinate[label=left:{$y$}] (a) at ($(one)!.4!(two)$);
\end{tikzpicture}
\raisebox{0.5cm}{$\rightarrow$}
\begin{tikzpicture}[scale=0.55]
\node (zero)[tn,label=0:{}]{}[grow=down]
	%[sibling distance=16mm,level distance=10mm]
    child[grow=231] {node [label=180:{}](one){}
     %[sibling distance=14mm,level distance=13mm]
            child {node [label=0:{}](two){}
            child {node [tn,label=0:{}](two2){}}
        }}
    child[grow=310] {node [label=0:{}](three){}
    child[grow=310] {node [tn,label=0:{}](four){}
    child[grow=55] {node [label=0:{}](five){}}}
        };
        \draw(zero)--(two2);
        \draw(zero)--(four);
       \node [tn,label=above:{$*$}] (1) at
(zero) {};
\coordinate[label=left:{$x$}] (x) at ($(zero)!.43!(two2)$);
\coordinate[label=right:{$x$}] (x) at ($(zero)!.43!(four)$);
\end{tikzpicture}\\
\hspace*{-0.4cm}$y\rightarrow x^2$
\end{center}
\end{minipage}\hspace{-1cm}
\begin{minipage}{0.51\linewidth}
\begin{center}
\hspace*{1.6cm}
\begin{tikzpicture}[scale=0.55]
\node (zero)[label=0:{}]{}[grow=down]
	%[sibling distance=16mm,level distance=10mm]
    child[grow=310] {node [tn, label=0:{}](one){}
    child[grow=310,dashed,line width=1pt] {node [tn, label=0:{}](two){}}
        };
     \coordinate[label=right:{$y$}] (x) at ($(one)!.4!(two)$);
\end{tikzpicture}
\hspace{1mm}
\raisebox{0.5cm}{$\rightarrow$}
\hspace{1mm}
\begin{tikzpicture}[scale=0.55]
\node (zero)[tn,label=0:{}]{}[grow=down]
	%[sibling distance=16mm,level distance=10mm]
    child[grow=231] {node [label=180:{}](one){}
     %[sibling distance=14mm,level distance=13mm]
            child {node [tn,label=0:{}](two){}
            child[grow=135] {node [label=0:{}](two2){}}
        }}
    child[grow=310] {node [label=0:{}](three){}
    child[grow=310] {node [label=0:{}](four){}
    child[grow=310] {node [tn,label=0:{}](five){}}}
        };
         \draw(zero)--(five);
         \draw(zero)--(two);
       \node [tn,label=above:{$*$}] (1) at
(zero) {};
\coordinate[label=right:{$x$}] (x) at ($(zero)!.45!(five)$);
\coordinate[label=left:{$x$}] (x) at ($(zero)!.45!(two)$);
\end{tikzpicture}\\
\hspace*{-0.2cm}$y\rightarrow x^2$
\end{center}
\end{minipage}
\caption{Insertion at a position labeled by $y$.}
\label{fy}
\end{figure}

If the position $i$ is labeled by $y$,
no matter whether it is in an
 up run or a down run, the change of weights is always
 consistent with the substitution rule $y\rightarrow x^2$,
 as illustrated in Figure \ref{fy}.

We are left with three cases regarding the last two elements of $\sigma$.
Keep in mind that if $\sigma$ ends with an up run, the last two labels
must be $ax$.

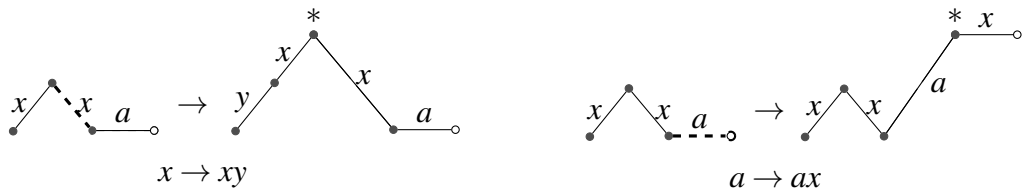
\begin{figure}[!ht]
\begin{minipage}{0.51\linewidth}
\begin{center}
\hspace*{1.3cm}
\begin{tikzpicture}[scale=0.55]
\node (zero)[tn,label=0:{}]{}[grow=down]
	%[sibling distance=16mm,level distance=10mm]
    child[grow=231] {node [tn, label=0:{}](one){}}
    child[grow=310,dashed] {node [tn, label=0:{}](two){}
    child[grow=0] {node [tn1, label=0:{}](three){}}
        };
    \coordinate[label=left:{$x$}] (a) at ($(zero)!.47!(one)$);
     \coordinate[label=right:{$x$}] (x) at ($(zero)!.47!(two)$);
          \coordinate[label=above:{$a$}] (x) at ($(two)!.5!(three)$);
         \node[tn1] (B) at (three){};
         \draw(two)--(three);
         \draw[dashed, line width=1.2pt](zero)--(two);
\end{tikzpicture}\hspace{0.3mm}
\raisebox{0.3cm}{$\rightarrow$}
\hspace{0.3mm}
\begin{tikzpicture}[scale=0.55]
\node (zero)[tn,label=0:{}]{}[grow=down]
	%[sibling distance=16mm,level distance=10mm]
    child[grow=231] {node [tn, label=180:{}](one){}
     %[sibling distance=14mm,level distance=13mm]
            child {node [tn,label=0:{}](two){}
        }}
    child[grow=310] {node [label=0:{}](three){}
    child[grow=310] {node [tn,label=0:{}](four){}
    child[grow=0] {node [tn1,label=0:{}](five){}}
    }
        };
         \draw(zero)--(four);
       \node [tn,label=above:{$*$}] (1) at
(zero) {};
\coordinate[label=left:{$x$}] (x) at ($(zero)!.4!(one)$);
    \coordinate[label=left:{$y$}] (y) at ($(one)!.4!(two)$);    \coordinate[label=right:{$x$}] (xx) at ($(zero)!.45!(four)$);
        \coordinate[label=above:{$a$}] (y) at ($(four)!.5!(five)$);
\end{tikzpicture}\\
\hspace*{0.5cm}$x\rightarrow xy$
\end{center}
\end{minipage}\hspace{-0.5cm}
\begin{minipage}{0.51\linewidth}
\begin{center}
\hspace*{1.4cm}
\begin{tikzpicture}[scale=0.55]
\node (zero)[tn,label=0:{}]{}[grow=down]
	%[sibling distance=16mm,level distance=10mm]
    child[grow=231] {node [tn, label=0:{}](one){}}
    child[grow=310] {node [tn, label=0:{}](two){}
    child[grow=0,dashed,line width=1.2pt] {node [tn1, label=0:{}](three){}}
        };
    \coordinate[label=left:{$x$}] (a) at ($(zero)!.5!(one)$);
     \coordinate[label=right:{$x$}] (x) at ($(zero)!.5!(two)$);
          \coordinate[label=above:{$a$}] (x) at ($(two)!.5!(three)$);
         \node[tn1] (B) at (three){};
\end{tikzpicture}
\raisebox{0.3cm}{$\rightarrow$}
\begin{tikzpicture}[scale=0.55]
\node (zero)[tn,label=0:{}]{}[grow=down]
	%[sibling distance=16mm,level distance=10mm]
    child[grow=231] {node [tn, label=180:{}](one){}
     }
    child[grow=310] {node [tn,label=0:{}](two){}
    child[grow=55] {node [label=0:{$a$}](three){}
    child[grow=55] {node [tn,label=0:{}](four){}
    child[grow=0] {node [tn1,label=0:{}](five){}}
       } }};
         \draw(two)--(four);                \node [tn,label=above:{$*$}] (1) at
(four) {};
\coordinate[label=left:{$x$}] (x) at ($(zero)!.4!(one)$);
    \coordinate[label=right:{$x$}] (y) at ($(zero)!.4!(two)$);
        \coordinate[label=above:{$x$}] (x2) at ($(four)!.5!(five)$);

\end{tikzpicture}\\
\hspace*{0.7cm}$a\rightarrow ax$
\end{center}
\end{minipage}
\caption{The labeling of $\sigma$ ends with $xa$.}
\label{fxa}
\end{figure}

\noindent
Case 1. $\sigma$ ends with a down run and the last two labels
are $x$ and $a$. As discussed before, the insertion at the position labeled by $x$
is in accordance  with the rule $x\rightarrow xy$ and the insertion at
the position of $a$ is reflected by the rule $a\rightarrow ax$, see
Figure \ref{fxa}.

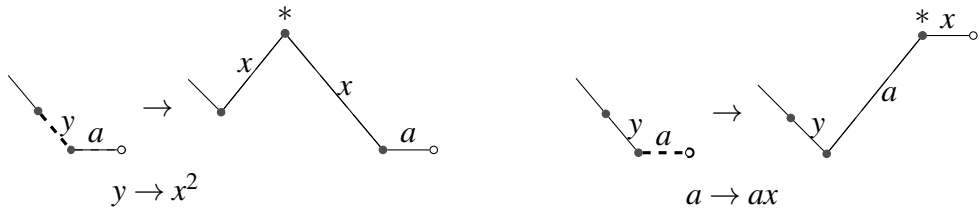
\begin{figure}[!ht]
\begin{minipage}{0.51\linewidth}
\begin{center}
\hspace*{1.3cm}
\begin{tikzpicture}[scale=0.45]
\node (zero)[label=0:{}]{}[grow=down]
	%[sibling distance=16mm,level distance=10mm]
    child[grow=310] {node [tn, label=0:{}](one){}
    child[dashed] {node [tn,label=180:{}](two){}
    child[grow=0] {node [tn1,label=0:{}](three){}}}
        };
     \coordinate[label=right:{$y$}] (x) at ($(one)!.4!(two)$);
          \coordinate[label=above:{$a$}] (a) at ($(two)!.5!(three)$);
          \node[tn1] at (three) {};
                    \draw(two)--(three);
          \draw[dashed, line width=1.2pt](one)--(two);
\end{tikzpicture}
\raisebox{0.5cm}{$\rightarrow$}
\begin{tikzpicture}[scale=0.45]
\node (zero)[tn,label=0:{}]{}[grow=down]
	%[sibling distance=16mm,level distance=10mm]
    child[grow=231] {node [label=180:{}](one){}
     %[sibling distance=14mm,level distance=13mm]
            child {node [tn,label=0:{}](two){}
            child[grow=135] {node [label=0:{}](two2){}}
        }}
    child[grow=310] {node [label=0:{}](three){}
    child[grow=310] {node [label=0:{}](four){}
    child[grow=310] {node [tn,label=0:{}](five){}
    child[grow=0] {node [tn1,label=0:{}](six){}}
    }}
        };
         \draw(zero)--(five);
         \draw(zero)--(two);
       \node [tn,label=above:{$*$}] (1) at
(zero) {};
\coordinate[label=left:{$x$}] (x) at ($(zero)!.4!(two)$);
\coordinate[label=right:{$x$}] (x) at ($(zero)!.47!(five)$);
\coordinate[label=above:{$a$}] (a) at ($(five)!.5!(six)$);
\end{tikzpicture}\\
\hspace*{-0.4cm}$y\rightarrow x^2$
\end{center}
\end{minipage}\hspace{-0.8cm}
\begin{minipage}{0.51\linewidth}
\begin{center}
\hspace*{1.6cm}
\begin{tikzpicture}[scale=0.45]
\node (zero)[label=0:{}]{}[grow=down]
	%[sibling distance=16mm,level distance=10mm]
    child[grow=310] {node [tn, label=0:{}](one){}
    child[grow=310] {node [tn, label=0:{}](two){}
     child[grow=0,dashed,line width=1.2pt] {node [tn1, label=0:{}](three){}}
    }
        };
     \coordinate[label=right:{$y$}] (x) at ($(one)!.4!(two)$);
     \coordinate[label=above:{$a$}] (a) at ($(two)!.5!(three)$);
     \node[tn1] (B) at (three){};
\end{tikzpicture}
\raisebox{0.5cm}{$\rightarrow$}
\begin{tikzpicture}[scale=0.45]
\node (zero)[tn,label=0:{}]{}[grow=down]
	%[sibling distance=16mm,level distance=10mm]
    child[grow=231] {node [label=180:{}](one){}
     %[sibling distance=14mm,level distance=13mm]
            child {node [label=0:{}](two){}
            child {node [tn,label=0:{}](three){}
            child[grow=135] {node [tn,label=0:{}](four){}
            child[grow=135] {node [label=0:{}](five){}}
            }
            }
        }}
    child[grow=0] {node [tn1, label=0:{}](six){}
        };
         \draw(zero)--(three);
       \node [tn,label=above:{$*$}] (1) at
(zero) {};
\coordinate[label=right:{$a$}] (x) at ($(zero)!.5!(three)$);
\coordinate[label=right:{$y$}] (y) at ($(three)!.7!(four)$);
\coordinate[label=above:{$x$}] (x) at ($(zero)!.5!(six)$);
\end{tikzpicture}\\
\hspace*{0.5cm}$a\rightarrow ax$
\end{center}
\end{minipage}
\caption{The labeling of $\sigma$ ends with $ya$.}
\label{fya}
\end{figure}

\noindent
Case 2. The labeling of $\sigma$ ends with $ya$.
As shown in Figure \ref{fya}, the changes of weights
caused by the insertions are coded by
the corresponding substitutions.

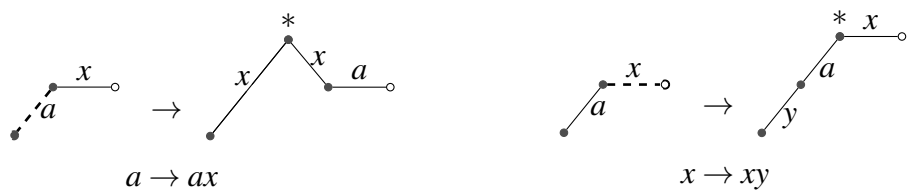
\begin{figure}[!ht]
\begin{minipage}{0.51\linewidth}
\begin{center}
\hspace*{1.3cm}
\begin{tikzpicture}[scale=0.55]
\node (zero)[tn,label=0:{}]{}[grow=down]
	%[sibling distance=16mm,level distance=10mm]
    child[grow=231,dashed,line width=1pt] {node [tn, label=0:{}](one){}}
    child[grow=0] {node [tn1, label=0:{}](two){}
        };
    \coordinate[label=right:{$a$}] (a) at ($(zero)!.5!(one)$);
     \coordinate[label=above:{$x$}] (x) at ($(zero)!.5!(two)$);
\end{tikzpicture}
\hspace{1mm}
\raisebox{0.3cm}{$\rightarrow$}
\hspace{1mm}
\begin{tikzpicture}[scale=0.55]
\node (zero)[tn,label=0:{}]{}[grow=down]
	%[sibling distance=16mm,level distance=10mm]
    child[grow=231] {node [label=180:{}](one){}
     %[sibling distance=14mm,level distance=13mm]
            child {node [tn,label=0:{}](two){}
        }}
    child[grow=310] {node [tn,label=0:{}](three){}
    child[grow=0] {node [tn1,label=0:{}](four){}}
        };
        \draw(zero)--(two);
       \node [tn,label=above:{$*$}] (1) at
(zero) {};
\coordinate[label=left:{$x$}] (x) at ($(zero)!.4!(two)$);
\coordinate[label=right:{$x$}] (x) at ($(zero)!.4!(three)$);
    \coordinate[label=above:{$a$}] (a) at ($(three)!.5!(four)$);
\end{tikzpicture}\\
\hspace*{0.5cm}$a\rightarrow ax$
\end{center}
\end{minipage}\hspace{-1cm}
\begin{minipage}{0.51\linewidth}
\begin{center}
\hspace*{1.4cm}
\begin{tikzpicture}[scale=0.55]
\node (zero)[tn,label=0:{}]{}[grow=down]
	%[sibling distance=16mm,level distance=10mm]
    child[grow=231] {node [tn, label=0:{}](one){}}
    child[grow=0,dashed,line width=1pt] {node [tn1, label=0:{}](two){}
        };
    \coordinate[label=right:{$a$}] (a) at ($(zero)!.5!(one)$);
     \coordinate[label=above:{$x$}] (x) at ($(zero)!.5!(two)$);
     \node[tn1] (B) at (two){};
\end{tikzpicture}
\hspace{1mm}
\raisebox{0.3cm}{$\rightarrow$}
\hspace{1mm}
\begin{tikzpicture}[scale=0.55]
\node (zero)[tn,label=0:{}]{}[grow=down]
	%[sibling distance=16mm,level distance=10mm]
    child[grow=231] {node [tn,label=180:{}](one){}
     %[sibling distance=14mm,level distance=13mm]
            child {node [tn,label=0:{}](two){}
        }}
    child[grow=0] {node [tn1,label=0:{}](three){}
        };
       \node [tn,label=above:{$*$}] (1) at
(zero) {};
\coordinate[label=right:{$a$}] (x) at ($(zero)!.7!(one)$);
    \coordinate[label=right:{$y$}] (y) at ($(one)!.7!(two)$);
        \coordinate[label=above:{$x$}] (a) at ($(zero)!.5!(three)$);
\end{tikzpicture}\\
\hspace*{1.2cm}$x\rightarrow xy$
\end{center}
\end{minipage}
\caption{The labeling of $\sigma$ ends with $ax$.}
\label{fa}
\end{figure}

\noindent
 Case 3. $\sigma$ ends with an up run, that is,
the last two labels are $a$ and $x$. As depicted in Figure \ref{fa},
in either case, the change of weights is in compliance with the grammar,
namely, the rules $a\rightarrow ax$ and $x \rightarrow xy$.

Summing up all the cases completes the proof. \qed

Once we have the above up-down labeling of a permutation,
we see that   $D^n(a)=a \Lambda_n(x,y)$ for all $n$.  The up-down
labeling also gives rise to a labeling  for
alternating runs. Let $n\geq 1$, and let $\sigma$  be a permutation
of $[n]$. Assume that $\sigma=\pi 1 \tau$.
Then the alternating labeling of $\sigma$ consists of the
up-down labelings of $\pi'$ and $\tau$, where $\pi'$ is the
reverse of $\pi$.
For example, the alternating labeling
of the permutation $7\, 3\, 5\,  8\,  6\,  1\,
4\,  9\,  2$ is given by
\[ (  x \  7\ a \   3\ y \ 5\ x\  8\ x \  6\  y \    1) \;
    (1 \ y\   4\ x \    9 \  x \  2 \  a).  \]
However, we shall not dwell further in this respect.

\section{Grammar assisted bijections}

In this section, we present a bijection
between permutations and increasing trees
that maps the number of up-down runs to the
number of even degree nonroot vertices of the
corresponding increasing tree. While the existence
of such a bijection is assured by
 the two grammatical labelings of the same grammar, a direct construction is less
obvious. We shall work out an explicit correspondence.

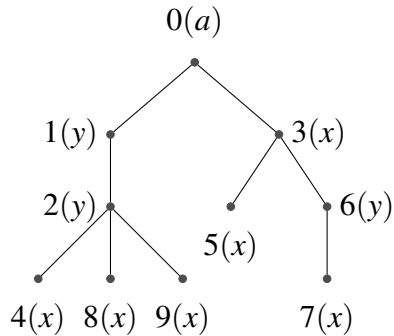
\begin{figure}[!ht]
\centering
\begin{tikzpicture}[scale=0.8]
\node (zero)[tn,label=90:{$0(a)$}]{}[grow=down]
	[sibling distance=28mm,level distance=12mm]
    child {node [tn, label=180:{$1(y)$}](one){}
    child {node [tn, label=180:{$2(y)$}](two){}
    [sibling distance=12mm,level distance=12mm]
    child {node [tn, label=270:{$4(x)$}](three){}}
    child {node [tn, label=270:{$8(x)$}](four){}}
    child {node [tn, label=270:{$9(x)$}](five){}}
    }
    }
    child {node [tn, label=0:{$3(x)$}](six){}
    [sibling distance=16mm,level distance=12mm]
        child {node [tn, label=270:{$5(x)$}](eight){}
    }
    child {node [tn, label=0:{$6(y)$}](seven){}
    child {node [tn, label=270:{$7(x)$}](nine){}}
    }
        };
\end{tikzpicture}
\caption{An increasing tree with the parity labeling.}
\label{fit}
\end{figure}

We begin with a labeling scheme for increasing trees, which
we call the parity labeling. Let $n\geq 1$, and let $T$ be an
increasing tree on $\{0, 1, \ldots, n\}$. The degree of a vertex
of $T$ is understood to be the number of its children. Then
the root of $T$ is labeled by $a$. If $v$ is not the root, then
it is labeled by $x$ if it is of even degree, and it is labeled
by $y$ if it is of odd degree. For example, an increasing tree
along with the parity labeling is shown in Figure \ref{fit}.

The grammar implies the following  bijection.

\begin{thm} \label{thm-bijection}
For $n\geq 1$, there is a bijection $\phi$ between the
set of permutations $\sigma$  of $[n]$ and the set
of increasing trees $T$ on $\{0,1, \ldots, n\}$ such that
the number of up-down runs of $\sigma$ equals the number of
even degree nonroot vertices of $\phi(\sigma)$.
\end{thm}

\begin{figure}[!ht]
\begin{center}
\raisebox{1.2cm}{$1\ 2\  \rightarrow\  0\  y\  1\  a\  2\  x \ \rightarrow$\quad }
\begin{tikzpicture}[scale=0.55]
\node (zero)[tn,label=0:{$0(a)$}]{}[grow=down]
    child {node [tn, label=0:{$1(y)$}](one){}
    child {node [tn, label=0:{$2(x)$}](two){}
    }
    };
\end{tikzpicture}
\hspace*{1.2cm}\raisebox{1.2cm}{$2\ 1\  \rightarrow\  0\  x\  2\  x\  1\  a \  \rightarrow$\ }
\begin{tikzpicture}[scale=0.55]
\node (zero)[tn,label=0:{$0(a)$}]{}[grow=down]
    child[grow=235] {node [tn, label=270:{$1(x)$}](one){}}
    child[grow=305] {node [tn, label=270:{$2(x)$}](two){}
    };
\end{tikzpicture}
\end{center}
\caption{The correspondence for $n=2$.}
\label{fn2}
\end{figure}
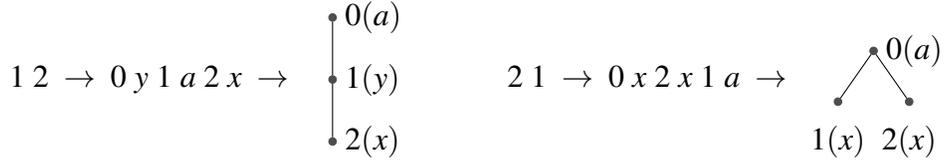

\noindent
{\it Proof. }
We proceed to describe the construction of the
map $\phi$.
For $n=1$, there is nothing to be said. For $n=2$, the correspondence is
unique subject to the weight preserving requirement, and it is shown in
Figure \ref{fn2}.

We now assume that $n\geq 3$. Let $\sigma$ be
a permutation of $[n]$, and let $\sigma^{(i)}$ be
the permutation obtained from $\sigma$ by removing
all the elements that are greater than $i$.
In order to construct an
increasing tree $T$ on $\{0, 1, \ldots, n\}$
from $\sigma$  with the weight preserving property,
we devise a procedure to generate a sequence
of increasing trees $T^{(2)}$, $T^{(3)}$, $\ldots$, $T^{(n)}$,
where $T^{(2)}$ is the increasing tree corresponding to
$\sigma^{(2)}$. As will be seen, for any $2\leq i \leq n$,
$\sigma^{(i)}$ and $T^{(i)}$ not only have the same weight,
but also share the same labeling that will be made clear later in the
course of the construction. We say that the
labeling of $\sigma^{(i)}$ is coherent with the labeling of
$T^{(i)}$ provided that the following conditions are satisfied. Set $\pi=\sigma^{(i)}$.
      \begin{itemize}
      \item If
      $\pi_{i-1} > \pi_i$, then for $1\leq k \leq i$, the position
      $k$ of $\pi$ has the same label
      as the vertex $\pi_k$ in $T^{(i)}$ and the position $i+1$ has
      the same label $a$ as the root of $T^{(i)}$.
      \item
       If
      $\pi_{i-1} < \pi_i$, then for $1\leq k \leq i-1$, the position
      $k$ of $\pi$ has the same label
      as the vertex $\pi_k$ in $T^{(i)}$,  the position $i$ has
      the same label $a$ as the root of $T^{(i)}$ and the position
      $i+1$ of $\pi$ has the same label $x$ as the vertex $\pi_{i}$ in
      $T^{(i)}$.
\end{itemize}

Let us now describe the construction of $T^{(i+1)}$ from
$T^{(i)}$.  Let $\pi=\sigma^{(i )}$.
Assume that $\sigma^{(i+1)}$ is obtained
from $\pi$ by inserting $i+1$ at
position $k$, where $1\leq k \leq i+1$.
We need to distinguish two main cases.
First, we assume that $\pi$ ends with a down run, that is,
$\pi$ has an even number of up-down runs, and so
the number of $x$ labels of $\pi$ is even.
Then the increasing tree $T^{(i+1)}$ is generated via
the following procedure.

\begin{enumerate}
\item Assume that
     the position $k$ of $\pi$ is labeled by $x$ and $k$ is
      on the rise, then
       the next position must be a down step labeled by $x$, which
       implies that  $k\leq i-1$. We call the position $k+1$ the dual
       position of $k$ and adjoin the vertex $i+1$ to the vertex $\pi_{k+1}$
      of $T^{(i)}$ as a child to obtain $T^{(i+1)}$ for which the
      label of $\pi_{k+1}$ changes from $x$ to $y$, see
      Figure \ref{fgx1}, where the square vertex
      signifies where to attach the vertex $i+1$ to
      $T^{(i)}$. Intuitively, for a position $k$
      labeled by $x$ that is on the rise, we need to look at
      the next position for corresponding operation on
      the increasing tree. Fortunately, one finds that
      the assumption that the position $k+1$ is labeled by
      $x$ is  the very property required for
      the procedure to work.
      Once it has been noticed, it is not hard to see that
      it is valid all along.

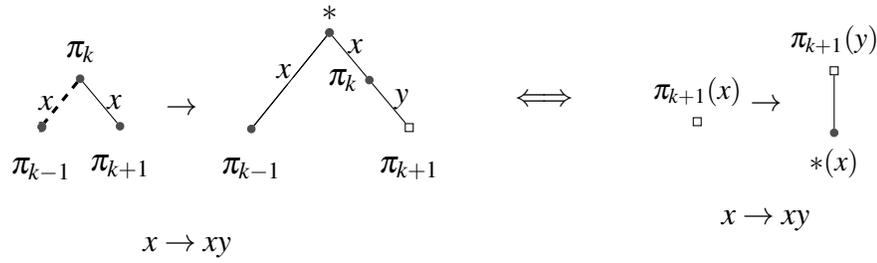
\begin{figure}[!ht]
\begin{minipage}{0.51\linewidth}
\begin{center}
\hspace*{0.3cm}
\begin{tikzpicture}[scale=0.55]
\node (zero)[tn,label=90:{$\pi_{k}$}]{}[grow=down]
	%[sibling distance=16mm,level distance=10mm]
    child[grow=231,dashed, line width=1.2pt] {node [tn, label=270:{$\pi_{k-1}$}](one){}}
    child[grow=310] {node [tn, label=270:{$\pi_{k+1}$}](two){}
        };
    \coordinate[label=left:{$x$}] (a) at ($(zero)!.5!(one)$);
     \coordinate[label=right:{$x$}] (x) at ($(zero)!.5!(two)$);
\end{tikzpicture}\hspace{0.3mm}
\raisebox{1.2cm}{$\rightarrow$}
\hspace{0.3mm}
\begin{tikzpicture}[scale=0.55]
\node (zero)[tn,label=0:{}]{}[grow=down]
	%[sibling distance=16mm,level distance=10mm]
    child[grow=231] {node [label=180:{}](one){}
     %[sibling distance=14mm,level distance=13mm]
            child {node [tn,label=270:{$\pi_{k-1}$}](two){}
        }}
    child[grow=310] {node [tn, label=180:{$\pi_{k}$}](three){}
    child[grow=310] {node [tn3,label=270:{$\pi_{k+1}$}](four){}}
        };
         \draw(zero)--(two);
       \node [tn,label=above:{$*$}] (1) at
(zero) {};
\coordinate[label=left:{$x$}] (x) at ($(zero)!.4!(two)$);
    \coordinate[label=right:{$y$}] (y) at ($(three)!.4!(four)$);
        \coordinate[label=right:{$x$}] (xr) at ($(zero)!.3!(three)$);
\end{tikzpicture}\\
\hspace*{-6mm}$x\rightarrow xy$
\end{center}
\end{minipage}
\raisebox{4mm}{$\Longleftrightarrow$}\hspace*{-15mm}
\begin{minipage}{0.51\linewidth}
\begin{center}
%\hspace*{1.4cm}
{\large $\substack{\pi_{k+1}(x)\\ \begin{tikzpicture}[scale=0.55]
       \node [tn3,label=above:{}] (1) at
(zero) {};
	%[sibling distance=16mm,level distance=10mm]
\end{tikzpicture}}$}
$\rightarrow$
{\large $\substack{\pi_{k+1}(y)\\
\begin{tikzpicture}[scale=0.55]
\node (zero)[tn3,label=above:{}]{}[grow=down]
	%[sibling distance=16mm,level distance=10mm]
    child[grow=270] {node [tn, label=270:{}](one){}};
\end{tikzpicture}\\
*(x)}$}\\
\vskip 4mm
$x\rightarrow xy$
\end{center}
\end{minipage}
\caption{$x$ is on the rise.}
\label{fgx1}
\end{figure}

\item Assume that
       the position $k$ of $\pi$ is labeled by $x$ and $k$ is
      on the fall. We call the position $k-1$ the dual
      position of $k$, and adjoin the vertex $i+1$ to the vertex $\pi_{k-1}$
      of $T^{(i )}$ as a child to obtain $T^{(i+1)}$, see Figure
      \ref{fgx2}, where and throughout the proof the symbol $*$ stands for the element $i+1$ to be added into $\sigma^{(i)}$ and
      $T^{(i)}$. As before, there is no danger to assume that
      the vertex $\pi_{k-1}$ is labeled by $x$.

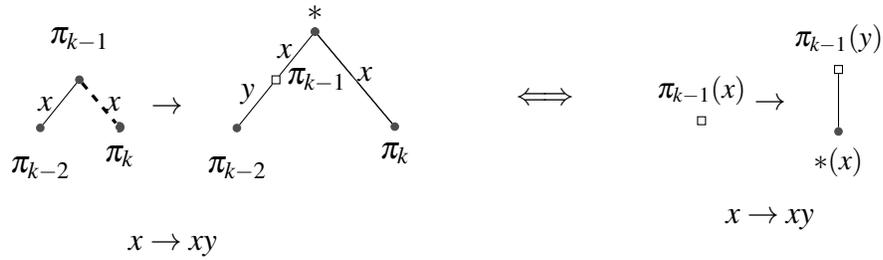
\begin{figure}[!ht]
\begin{minipage}{0.6\linewidth}
\begin{center}
\hspace*{0.3cm}
\begin{tikzpicture}[scale=0.55]
\node (zero)[tn,label=90:{$\pi_{k-1}$}]{}[grow=down]
	%[sibling distance=16mm,level distance=10mm]
    child[grow=231] {node [tn, label=270:{$\pi_{k-2}$}](one){}}
    child[grow=310,dashed, line width=1.2pt] {node [tn, label=270:{$\pi_{k}$}](two){}
        };
    \coordinate[label=left:{$x$}] (a) at ($(zero)!.5!(one)$);
     \coordinate[label=right:{$x$}] (x) at ($(zero)!.5!(two)$);
\end{tikzpicture}\hspace{0.3mm}
\raisebox{1.2cm}{$\rightarrow$}
\hspace{0.3mm}
\begin{tikzpicture}[scale=0.55]
\node (zero)[tn,label=0:{}]{}[grow=down]
	%[sibling distance=16mm,level distance=10mm]
    child[grow=231] {node [tn3, label=0:{$\pi_{k-1}$}](one){}
     %[sibling distance=14mm,level distance=13mm]
            child {node [tn,label=270:{$\pi_{k-2}$}](two){}
        }}
    child[grow=310] {node [label=0:{}](three){}
    child[grow=310] {node [tn,label=270:{$\pi_{k}$}](four){}}
        };
         \draw(zero)--(four);
       \node [tn,label=above:{$*$}] (1) at
(zero) {};
\coordinate[label=left:{$x$}] (x) at ($(zero)!.4!(one)$);
    \coordinate[label=left:{$y$}] (y) at ($(one)!.3!(two)$);
        \coordinate[label=right:{$x$}] (xr) at ($(zero)!.45!(four)$);
\end{tikzpicture}\\
\hspace*{-6mm}$x\rightarrow xy$
\end{center}
\end{minipage}
\hspace*{-6mm}
\raisebox{4mm}{$\Longleftrightarrow$}\hspace*{-6mm}
\begin{minipage}{0.4\linewidth}
\begin{center}
%\hspace*{1.4cm}
{\large $\substack{\pi_{k-1}(x)\\ \begin{tikzpicture}[scale=0.55]
       \node [tn3,label=above:{}] (1) at
(zero) {};
	%[sibling distance=16mm,level distance=10mm]
\end{tikzpicture}}$}
$\rightarrow$
{\large $\substack{\pi_{k-1}(y)\\
\begin{tikzpicture}[scale=0.55]
\node (zero)[tn3,label=above:{}]{}[grow=down]
	%[sibling distance=16mm,level distance=10mm]
    child[grow=270] {node [tn, label=270:{}](one){}};
\end{tikzpicture}\\
*(x)}$}\\
\vskip 4mm
$x\rightarrow xy$
\end{center}
\end{minipage}
\caption{$x$ is on the fall.}
\label{fgx2}
\end{figure}

\item Assume that the position $k$ of $\pi$ is labeled by $y$.
Then adjoin $i+1$ to the vertex $\pi_k$ in $T^{(i)}$ as
a child, see Figures \ref{fgy1} and \ref{fgy2}.

\begin{figure}[!ht]
\begin{minipage}{0.6\linewidth}
\begin{center}
\hspace*{0.3cm}
\begin{tikzpicture}[scale=0.5]
\node (zero)[tn,label=0:{$\pi_{k+1}$}]{}[grow=down]
	%[sibling distance=16mm,level distance=10mm]
    child[grow=231] {node [tn3, label=0:{$\pi_{k}$}](one){}
    child[grow=231,dashed, line width=1.2pt] {node [tn, label=270:{$\pi_{k-1}$}](two){}}
        };
     \coordinate[label=left:{$y$}] (x) at ($(one)!.4!(two)$);
\end{tikzpicture}\hspace{0.3mm}
\raisebox{1.6cm}{$\rightarrow$}
\hspace{0.3mm}
\begin{tikzpicture}[scale=0.50]
\node (zero)[tn,label=0:{}]{}[grow=down]
	%[sibling distance=16mm,level distance=10mm]
    child[grow=231] {node [label=0:{}](one){}
     %[sibling distance=14mm,level distance=13mm]
            child {node [label=270:{}](two){}
             child {node [tn,label=270:{$\pi_{k-1}$}](three){}}
        }}
    child[grow=310] {node [label=0:{}](four){}
    child[grow=310] {node [tn3,label=270:{$\pi_{k}$}](five){}
    child[grow=55] {node [tn,label=0:{$\pi_{k+1}$}](six){}}
    }
        };
         \draw(zero)--(three);
         \draw(zero)--(five);
       \node [tn,label=above:{$*$}] (1) at
(zero) {};
\coordinate[label=left:{$x$}] (x) at ($(zero)!.45!(three)$);
    \coordinate[label=right:{$x$}] (xy) at ($(zero)!.45!(five)$);
\end{tikzpicture}\\
\hspace*{-13mm}$y\rightarrow x^2$
\end{center}
\end{minipage}
\raisebox{4mm}{$\Longleftrightarrow$}\hspace*{-10mm}
\begin{minipage}{0.4\linewidth}
\begin{center}
%\hspace*{1.4cm}
{\large $\substack{\pi_{k}(y)\\ \begin{tikzpicture}[scale=0.55]
       \node [tn3,label=above:{}] (1) at
(zero) {};
	%[sibling distance=16mm,level distance=10mm]
\end{tikzpicture}}$}
$\rightarrow$
{\large $\substack{\pi_{k}(x)\\
\begin{tikzpicture}[scale=0.55]
\node (zero)[tn3,label=above:{}]{}[grow=down]
	%[sibling distance=16mm,level distance=10mm]
    child[grow=270] {node [tn, label=270:{}](one){}};
\end{tikzpicture}\\
*(x)}$}\\
\vskip 4mm
$y\rightarrow x^2$
\end{center}
\end{minipage}
\caption{$y$ is on the rise.}
\label{fgy1}
\end{figure}
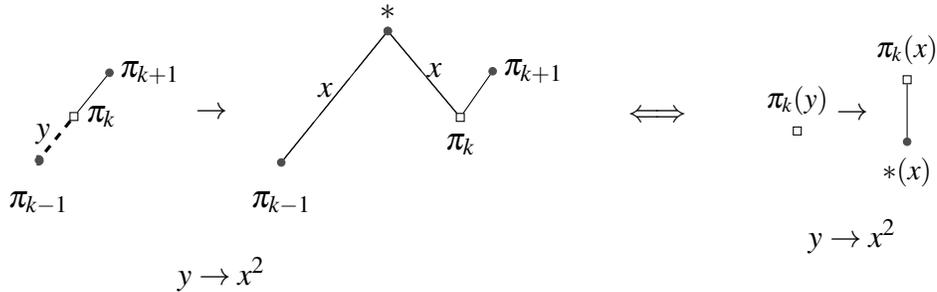

\begin{figure}[!ht]
\begin{minipage}{0.6\linewidth}
\begin{center}
\hspace*{0.3cm}
\begin{tikzpicture}[scale=0.55]
\node (zero)[tn,label=90:{$\pi_{k-2}$}]{}[grow=down]
	%[sibling distance=16mm,level distance=10mm]
    child[grow=310] {node [tn, label=180:{$\pi_{k-1}$}](one){}
    child[grow=310,dashed] {node [tn3, label=270:{$\pi_{k}$}](two){}}
        };
        \draw[dashed,line width=1.2pt](one)--(two);
    \node [tn3,label=right:{}] (1) at
(two) {};
     \coordinate[label=right:{$y$}] (x) at ($(one)!.4!(two)$);
\end{tikzpicture}\hspace{0.3mm}
\raisebox{1.6cm}{$\rightarrow$}
\hspace{0.3mm}
\begin{tikzpicture}[scale=0.55]
\node (zero)[tn,label=0:{}]{}[grow=down]
	%[sibling distance=16mm,level distance=10mm]
    child[grow=231] {node [label=0:{}](one){}
     %[sibling distance=14mm,level distance=13mm]
            child {node [tn,label=270:{$\pi_{k-1}$}](two){}
             child[grow=135] {node [tn,label=90:{$\pi_{k-2}$}](three){}}
        }}
    child[grow=310] {node [label=0:{}](four){}
    child[grow=310] {node [label=270:{}](five){}
    child[grow=310] {node [tn3,label=270:{$\pi_{k}$}](six){}}
    }
        };
         \draw(zero)--(two);
         \draw(zero)--(six);
       \node [tn,label=above:{$*$}] (1) at
(zero) {};
\coordinate[label=left:{$x$}] (x) at ($(zero)!.45!(two)$);
    \coordinate[label=right:{$x$}] (xy) at ($(zero)!.45!(six)$);
\end{tikzpicture}\\
$y\rightarrow x^2$\hspace*{18mm}
\end{center}
\end{minipage}
\hspace*{-2mm}\raisebox{4mm}{$\Longleftrightarrow$}\hspace*{-9mm}
\begin{minipage}{0.4\linewidth}
\begin{center}
%\hspace*{1.4cm}
{\large $\substack{\pi_{k}(y)\\ \begin{tikzpicture}[scale=0.55]
       \node [tn3,label=above:{}] (1) at
(zero) {};
	%[sibling distance=16mm,level distance=10mm]
\end{tikzpicture}}$}
$\rightarrow$
{\large $\substack{\pi_{k}(x)\\
\begin{tikzpicture}[scale=0.55]
\node (zero)[tn3,label=above:{}]{}[grow=down]
	%[sibling distance=16mm,level distance=10mm]
    child[grow=270] {node [tn, label=270:{}](one){}};
\end{tikzpicture}\\
*(x)
}$}\\
\vskip 4mm
$y\rightarrow x^2$
\end{center}
\end{minipage}
\caption{$y$ is on the fall.}
\label{fgy2}
\end{figure}
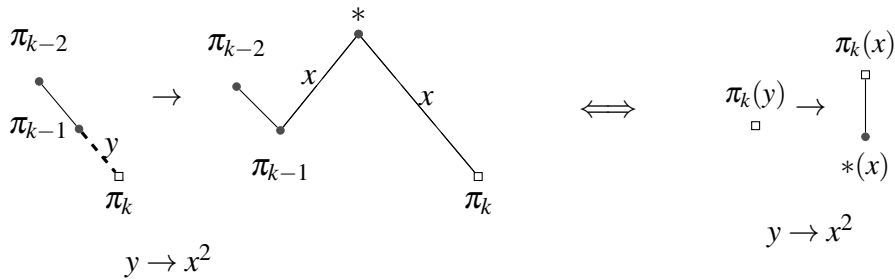
\item If $k=i+1$, then adjoin $i+1$ to the root of $T^{(i)}$.
In this case, the involved substitution is $a\rightarrow ax$.
\end{enumerate}

In any of the above cases,
it is observable that after the update of $T^{(i)}$,
the labeling of $\sigma^{(i+1)}$ is consistent with
 the labeling of $T^{(i+1)}$.
In other words, all the
assumptions in the above argument are well-grounded.

We are now left with the case that $\pi$ has an odd number
of up-down runs, that is, the labels of $\pi$
end with $ax$.   Here we encounter three possibilities.
\begin{enumerate}
  \item If $i+1$ is inserted at position $i$, namely,
  at the position before $\pi_{i}$, then adjoin the vertex $i+1$ to the root $0$.
 In this case, the vertex
 $i+1$ is labeled by $x$ and the root remains to be labeled by $a$, see Fig. 12.

\begin{figure}[!ht]
\begin{minipage}{0.6\linewidth}
\begin{center}
\hspace*{0.3cm}
\begin{tikzpicture}[scale=0.55]
\node (zero)[tn,label=90:{$\pi_{i}$}]{}[grow=down]
	%[sibling distance=16mm,level distance=10mm]
    child[grow=231,dashed, line width=1.2pt] {node [tn, label=270:{$\pi_{i-1}$}](one){}}
    child[grow=0] {node [tn3, label=90:{$0$}](two){}
        };
    \coordinate[label=left:{$a$}] (a) at ($(zero)!.5!(one)$);
     \coordinate[label=above:{$x$}] (x) at ($(zero)!.5!(two)$);
\end{tikzpicture}\hspace{0.3mm}
\raisebox{1.2cm}{$\rightarrow$}
\hspace{0.3mm}
\begin{tikzpicture}[scale=0.55]
\node (zero)[tn,label=0:{}]{}[grow=down]
	%[sibling distance=16mm,level distance=10mm]
    child[grow=231] {node [label=180:{}](one){}
     %[sibling distance=14mm,level distance=13mm]
            child {node [tn,label=270:{$\pi_{i-1}$}](two){}
        }}
    child[grow=310] {node [tn, label=270:{$\pi_{i}$}](three){}
    child[grow=0] {node [tn3,label=270:{$0$}](four){}}
        };
         \draw(zero)--(two);
       \node [tn,label=above:{$*$}] (1) at
(zero) {};
\coordinate[label=left:{$x$}] (x) at ($(zero)!.4!(two)$);
    \coordinate[label=above:{$a$}] (y) at ($(three)!.5!(four)$);
        \coordinate[label=right:{$x$}] (xr) at ($(zero)!.3!(three)$);
\end{tikzpicture}\\
\hspace*{-6mm}
$a\rightarrow ax$
\end{center}
\end{minipage}
\hspace*{-6mm}\raisebox{4mm}{$\Longleftrightarrow$}\hspace*{-10mm}
\begin{minipage}{0.4\linewidth}
\begin{center}
%\hspace*{1.4cm}
{\large $\substack{0(a)\\ \begin{tikzpicture}[scale=0.55]
       \node [tn3,label=above:{}] (1) at
(zero) {};
	%[sibling distance=16mm,level distance=10mm]
\end{tikzpicture}}$}
\hspace*{2mm}
$\rightarrow$\hspace*{2mm}
{\large $\substack{0(a)\\
\begin{tikzpicture}[scale=0.55]
\node (zero)[tn3,label=above:{}]{}[grow=down]
	%[sibling distance=16mm,level distance=10mm]
    child[grow=270] {node [tn, label=270:{}](one){}};
\end{tikzpicture}\\
*(x)}$}\\
\vskip 4mm
$a\rightarrow ax$
\end{center}
\end{minipage}
\caption{The case $k=i$ and $\pi_{i-1}< \pi_{i}$.}
\label{f21}
\end{figure}
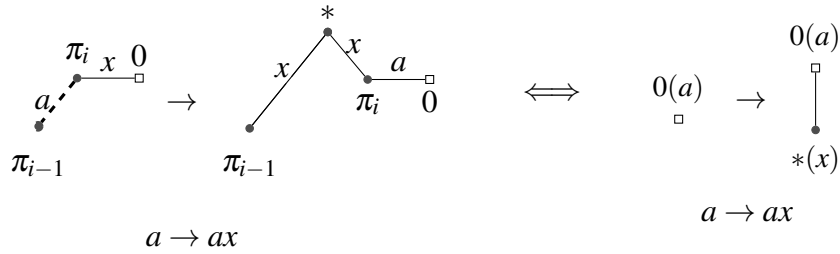

 \item
If $i+1$
 inserted into $\pi$ at position $i+1$, that is, $i+1$ is
  inserted at the end of $\pi$,
  then adjoin $i+1$ to the vertex $\pi_i$ of $T^{(i)}$. In this case,
 $i+1$ is labeled by $x$ in $T^{(i+1)}$ and the the
 label of $\pi_i$ is switched from $x$ to $y$, see Fig. 13.

\begin{figure}[!ht]
\begin{minipage}{0.6\linewidth}
\begin{center}
\hspace*{0.3cm}
\begin{tikzpicture}[scale=0.55]
\node (zero)[tn3,label=90:{$\pi_{i}$}]{}[grow=down]
	%[sibling distance=16mm,level distance=10mm]
    child[grow=231] {node [tn, label=270:{$\pi_{i-1}$}](one){}}
    child[grow=0,dashed] {node [tn2, label=90:{$0$}](two){}
        };
    \coordinate[label=left:{$a$}] (a) at ($(zero)!.4!(one)$);
     \coordinate[label=above:{$x$}] (x) at ($(zero)!.5!(two)$);
     \draw[dashed,line width=1.2pt](zero)--(two);
         \node [tn2,label=right:{}] (1) at
(two) {};
\end{tikzpicture}\hspace{0.3mm}
\raisebox{1.2cm}{$\rightarrow$}
\hspace{0.3mm}
\begin{tikzpicture}[scale=0.55]
\node (zero)[tn,label=0:{}]{}[grow=down]
	%[sibling distance=16mm,level distance=10mm]
    child[grow=231] {node [tn3,label=0:{$\pi_{i}$}](one){}
     %[sibling distance=14mm,level distance=13mm]
            child {node [tn,label=270:{$\pi_{i-1}$}](two){}
        }}
    child[grow=0] {node [tn2, label=90:{}](three){}
        };
         %\draw(zero)--(two);
       \node [tn,label=above:{$*$}] (1) at
(zero) {};
\coordinate[label=left:{$a$}] (x) at ($(zero)!.4!(one)$);
    \coordinate[label=left:{$y$}] (y) at ($(one)!.4!(two)$);
        \coordinate[label=above:{$x$}] (xr) at ($(zero)!.5!(three)$);
\end{tikzpicture}\\
%\hspace*{-3mm}
$x\rightarrow xy$
\end{center}
\end{minipage}
\hspace*{-10mm}\raisebox{4mm}{$\Longleftrightarrow$}\hspace*{-9mm}
\begin{minipage}{0.4\linewidth}
\begin{center}
%\hspace*{1.4cm}
{\large $\substack{\pi_i(x)\\ \begin{tikzpicture}[scale=0.55]
       \node [tn3,label=above:{}] (1) at
(zero) {};
	%[sibling distance=16mm,level distance=10mm]
\end{tikzpicture}}$}
\hspace*{2mm}
$\rightarrow$\hspace*{2mm}
{\large $\substack{\pi_i(y)\\
\begin{tikzpicture}[scale=0.55]
\node (zero)[tn3,label=above:{}]{}[grow=down]
	%[sibling distance=16mm,level distance=10mm]
    child[grow=270] {node [tn, label=270:{}](one){}};
\end{tikzpicture}\\
*(x)}$}\\
\vskip 4mm
\hspace*{3mm}
$x\rightarrow xy$
\end{center}
\end{minipage}
\caption{The case $k=i+1$ and $\pi_{i-1}<\pi_{i}$.}
\label{f22}
\end{figure}
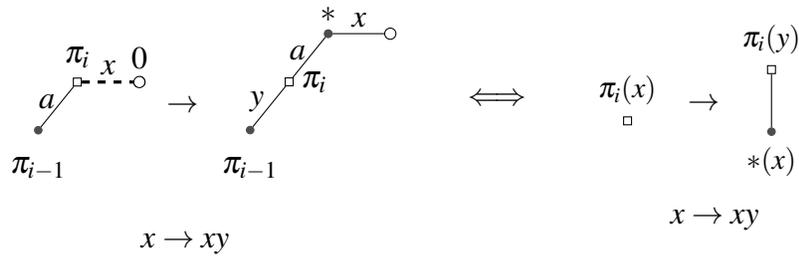

 \item If $i+1$ is inserted at position $k$ with $k\leq i-1$, then
  we may follow the above procedure for the case when $\pi$ has
  an even number of up-down runs to generate an increasing
  tree $T^{(i+1)}$.  It should be pointed out that there are no worries
  even when the position $k$ is on the
  rise and labeled by $y$, in which case $\pi_1\pi_2\cdots \pi_{k}$
  has an odd number of up-down runs,
  since under this circumstance the $x$ label still appear
  in pairs in $\pi$ before the position $k$.

\end{enumerate}

Up till now,  we have provided a procedure to build an
increasing tree $T$ from a permutation $\sigma$.
To see that the process is reversible, we may
extract from $T$  a sequence of increasing trees
$T^{(2)}$, $T^{(3)}$, $\ldots$, $T^{(n)}$, where
$T^{(i)}$ is obtained from $T$ by removing all the
vertices that are greater than $i$.

As
 $T^{(2)}$ is in one-to-one correspondence with a
 permutation $\sigma^{(2)}$, we may carry out the
 following steps by starting with $i=2$. Using
  $\sigma^{(i)}$ and $T^{(i+1)}$ along with the
  labeling assumption, one sees that $\sigma^{(i+1)}$
  can be retrieved, however the detailed reasoning
  is omitted. So we reach the conclusion that
  the procedure leads to a desired bijection. \qed

The operation of locating the vertex in an increasing tree
when an insertion into a permutation
takes place at a position labeled by $x$ is
called the reflection principle.
For example, for the increasing tree in Figure \ref{fit}, the
intermediate permutations are given in the table below, where
an underlined label indicates where the element $i+1$ is
inserted into $\sigma^{(i)}$.
Notice that $\sigma$ has six up-down runs, whereas
$T$ has six nonroot vertices of even degree.

\begin{center}
\renewcommand\arraystretch{1.33}
\begin{tabular}{|l|l|l|l|}
\hline
 $i$ \quad &  $\sigma^{(i)}$ with labeling & Weight  & Substitution\\ \hline
 2 &  $0\  y\  1\ \underline{a} \ 2\  x $ & $a x y$  &  $a\rightarrow ax$\\\hline
  3 & $0\  y \  1\ {x}\ 3\ \underline{x} \ 2 \  a$
  & $a x^2 y$ &  $x \rightarrow xy$ \\\hline
   4 &  $0\  y\ 1\  \underline{x}\ 4 \  {x} \  3\ y \ 2\ a$ & $a x^2 y^2$ & $x\rightarrow xy$ \\\hline
    5 &  $0\  y \  1\ x\  5\ x\ 4\ \underline{y}\ 3\ y\
     2\  a $  & $a x^2 y^3$ & $y\rightarrow x^2$ \\\hline
 6 & $0\  y\  1\ x \ 5\ x\ 4\ x\ 6\ \underline{x}\  3\  {y}\ 2\ a $     & $a x^4 y^2$  &  $x\rightarrow xy$\\\hline
  7 & $0\  y\  1\ x \ 5\ x\ 4\ y\ 6\ x\ 7\ x\  3\  \underline{y}\ 2\ a $  & $a x^4 y^3$  &$y \rightarrow x^2$ \\\hline
   8 &  $0\  y\  1\ {x}\ 5\ x\ 4\ y\ 6\ x\ 7\ x\ 3\ \underline{x}\ 8\ x\ 2\ a $    & $a x^6 y^2$ &
   $x\rightarrow xy$\\\hline
    9\; &    $0\  y\  1\ x\ 5 \ x \  4\ y\ 6\ x\ 7\ x\ 3\ x\ 9\ x\ 8\ y\ 2\ a $
       \quad \quad  & $a x^6 y^3$  &  \\\hline
\end{tabular}
\end{center}

It should be noted that the above bijection reduces to a
correspondence between alternating (down-up) permutations
of $[n]$ and
even increasing trees on $\{0, 1, \ldots, n\}$.
Recall that an alternating permutation $[n]$
is referred to a permutation $\sigma_1\sigma_2\cdots \sigma_n$
such that $\sigma_1>\sigma_2<\sigma_3> \cdots$, see
Stanley \cite{Stanley-2010}, whereas
an increasing tree is called even if every vertex possibly except
the root is of even degree. Notice that by the alternating condition
 some authors mean the up-down condition $\sigma_1< \sigma_2>\sigma_3< \cdots$ instead.
 Clearly, a permutation
 is alternating if and only if it has no $y$ labels
 in the up-down labeling, and down-up permutations are
 in one-to-one correspondence with up-down permutations via
 complementation.
 A bijection between up-down
permutations and even increasing trees has been given by
Kuznetsov, Pak and Postnikov \cite{KPP-1994}.

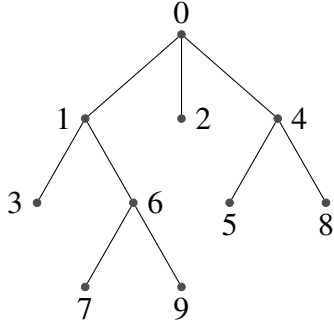
\begin{figure}[!ht]
\begin{center}
\begin{tikzpicture}[scale=0.8]
\node (zero)[tn,label=90:{$0$}]{}[grow=down]
[sibling distance=16mm,level distance=14mm]
    child {node [tn, label=180:{$1$}](one){}
    child {node [tn, label=180:{$3$}](three){}}
    child {node [tn, label=0:{$6$}](six){}
    child {node [tn, label=270:{$7$}](seven){}}
    child {node [tn, label=270:{$9$}](nine){}}
    }
    }
    child {node [tn, label=0:{$2$}](two){}}
    child {node [tn, label=0:{$4$}](four){}
    child {node [tn, label=270:{$5$}](five){}}
    child {node [tn, label=270:{$8$}](eight){}}
    };
\end{tikzpicture}
\end{center}
\caption{An even increasing tree.}
\label{f-p-i-e}
\end{figure}

When applied to a down-up permutation, the
procedure in the above
theorem produces an even increasing tree. For example, below is a down-up permutation together with
the up-down labeling:
\[ 0\ x \  3 \  x\   2\  x \  9\ x\  6\ x\  7 \ x\ 1\ x\ 8\ x\
4\  a\  5\  x.\]
The corresponding even increasing tree is displayed in Figure
\ref{f-p-i-e}.

We remark that a variant of the
construction in the above
proof gives a grammar assisted bijection linking
the number of left peaks of a permutation to the number
of  vertices of even degree  of an increasing tree, as given
in \cite{Chen-Fu-2017}.

\begin{thm}[\cite{Chen-Fu-2017}] \label{thm-L-I}
For $n \geq 1$ and $0\leq m \leq \lfloor n/2 \rfloor$, there
is a one-to-one correspondence between the set of
permutations of $[n]$ with $m$ left peaks and the
set of increasing trees on $\{0,1,\ldots, n\}$ with
$2m+1$ nonroot vertices of even degree.
\end{thm}

The grammar assisted bijection can be described as follows.
For $n\geq 2$, given a permutation $\sigma$ of $[n]$,
let $\sigma^{(2)}$, $\sigma^{(3)}$, $\ldots$, $\sigma^{(n)}=\sigma$
be the sequences of permutations such that $\sigma^{(i)}$ is
obtained from $\sigma$ by removing all the elements that are greater
than $i$. We wish to generate a sequence of
increasing trees $T^{(2)}$, $T^{(3)}$, $\ldots$, $T^{(n)}$ such
that $\sigma^{(i)}$ and $T^{(i)}$ have the same weight, as
clarified below.

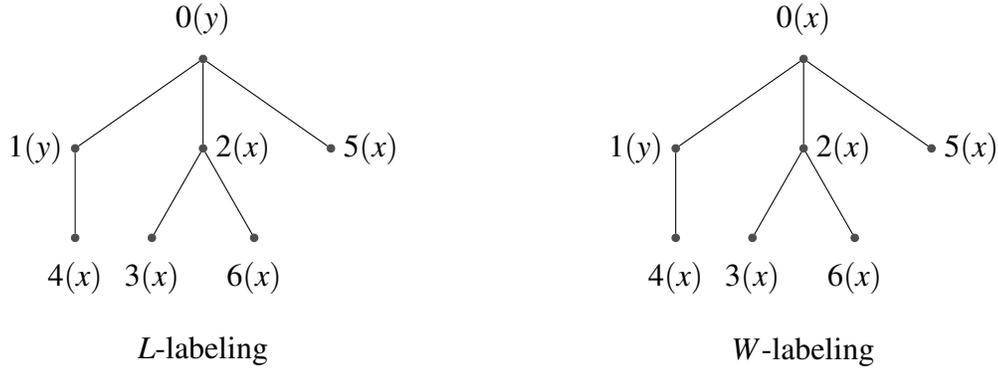
\begin{figure}[!ht]
\begin{minipage}{0.51\linewidth}
\begin{center}
\begin{tikzpicture}[scale=0.85]
\node (zero)[tn,label=90:{$0(y)$}]{}[grow=down]
[sibling distance=20mm,level distance=14mm]
    child {node [tn, label=180:{$1(y)$}](one){}
    child {node [tn, label=270:{$4(x)$}](four){}}
    }
    child {node [tn, label=0:{$2(x)$}](two){}
    [sibling distance=16mm,level distance=14mm]
    child {node [tn, label=270:{$3(x)$}](three){}}
    child {node [tn, label=270:{$6(x)$}](six){}}
    }
        child {node [tn, label=0:{$5(x)$}](five){}
    };
\end{tikzpicture}\\
$L$-labeling
\end{center}
\end{minipage}
\begin{minipage}{0.51\linewidth}
\begin{center}
\begin{tikzpicture}[scale=0.85]
\node (zero)[tn,label=90:{$0(x)$}]{}[grow=down]
[sibling distance=20mm,level distance=14mm]
    child {node [tn, label=180:{$1(y)$}](one){}
    child {node [tn, label=270:{$4(x)$}](four){}}
    }
    child {node [tn, label=0:{$2(x)$}](two){}
    [sibling distance=16mm,level distance=14mm]
    child {node [tn, label=270:{$3(x)$}](three){}}
    child {node [tn, label=270:{$6(x)$}](six){}}
    }
        child {node [tn, label=0:{$5(x)$}](five){}
    };
\end{tikzpicture}\\
$W$-labeling
\end{center}
\end{minipage}
\caption{Two labelings of an increasing tree.}
\label{FLW}
\end{figure}

Recall that the $L$-labeling of $\sigma$ is defined to
assign the label $x$ to the two positions that form a left peak
as well as to the
last position, and the rest of the positions are endowed with
 the label $y$.
For example, below is the $L$-labeling of a permutation of $[9]$:
\[ 0\  y\  1\ x\ 9 \ x\ 8\ y\ 3\ x\ 6\ x\  5\  y \ 4\ y\ 2\ y\ 7\  x. \]
On the other hand, for an increasing tree $T$ on $\{0,1,\ldots,n\}$,
a vertex is labeled by $x$ is of even degree, otherwise it is
labeled by $y$. We call this labeling scheme the $L$-labeling of $T$, see Figure
\ref{FLW}.

Suppose that
$\sigma^{(i+1)}$ is obtained from $\sigma^{(i)}$
by inserting $i+1$ at position $k$, where $1\leq k \leq i+1$.
Let $\pi=\sigma^{(i)}$.
 Observe that
the last label of $\sigma^{(i)}$ is always $x$ and the rest
of the $x$ labels always appear in pairs.
To produce $T^{(i+1)}$ from $T^{(i)}$, we stand by the following rules:
\begin{enumerate}
\item   Assume that $k\leq i-1$. Depending upon
whether the position $k$ is labeled by $x$ or $y$,
we proceed as in the case for up-down runs.

\item Assume that  $k=i$.  In this case,
the position $k$ is
labeled by $y$. If  $\pi_{i-1} < \pi_{i}$,
then adjoin $i+1$ to the root in $T^{(i)}$
to generate $T^{(i+1)}$; otherwise adjoin $i+1$ to
the vertex $\pi_i$.

\item Assume that $k=i+1$. If $\pi_{i-1}<\pi_i$, then
adjoin $i+1$ to the vertex $\pi_{i}$ of $T^{(i)}$ ; otherwise adjoin
$i+1$ to the root to
generate $T^{(i+1)}$.
\end{enumerate}

For instance, let us consider
the  permutation taken from the
 example for Theorem \ref{thm-bijection},  \[
\sigma=  1\ 5 \  4\ 6\   7\ 3\ 9\ 8\ 2,   \]
it can be checked that $\sigma$ corresponds to the same increasing tree as in
Figure \ref{fit} with the label $a$ of the root   replaced
by $x$. This is no coincidence, as will be seen.

Naturally, one may wonder
how to map the number of exterior   peaks
to a statistic of an increasing tree.
Here is a key observation about the $W$-labeling,
that is, the $x$  labels always
appear in pairs. This property allows us to
present the above grammar-assisted bijections in a much
more concise manner.

Imagine that an increasing tree with a single vertex $0$
is labeled by $y$. Then we comply with the parity rules
as before to produce an increasing tree
by successively adding vertices. The associated labeling
is called the $W$-labeling of an increasing tree, see Figure
\ref{FLW}.

\begin{thm}\label{thm-W-I}
For $n\geq 1$, there is a bijection mapping a permutation $\sigma$
of $[n]$ with $k$ exterior peaks to an increasing tree $T$ on $\{0, 1, \ldots, n\}$
with $j$ vertices of even degree such that $k = \lfloor {(j+1) / 2} \rfloor$.
\end{thm}

Note that for a permutation $\sigma$ with $k$ exterior peaks, the corresponding
increasing tree has $2k$ vertices labeled by $x$. This parity property becomes
transparent if we define the degree of a vertex of an increasing tree
as the number of adjacent vertices as if the tree is a graph
in the usual sense. Obviously
such a concern is  superficial,
and we had better stick to the usual terminology for rooted trees
especially
when we have a clear picture in mind.

We conclude  with a remark that
 all the aforementioned bijections
can be presented in a unified way  solely in terms of the
reflection principle. That is to say, the
bijection in
Theorem \ref{thm-W-I}
does the same job as the other grammar assisted bijections.
 For example,
the increasing tree in Figure \ref{FLW} invariably
corresponds to the permutation
$6\, 2\, 4\, 3\, 1\, 5$ for any of the
three labeling schemes.

We now return to the key step of the reflection operation.
Let $n \geq 1$, assume that $\sigma=\sigma_1\sigma_2\cdots \sigma_{n}$
is a permutation of $[n]$. For $1\leq k \leq n+1$, the position before
$\sigma_k$ is called position $k$, whereas the last position is
referred to position $0$. Now, the exterior peaks gather the corresponding
positions in pairs, so that we do not have to take special care of the last two elements
concerning their relative order.
 Moreover, we may index the positions of $\sigma$ by $0, 1, \ldots, n$,
and this yields a bijection that is seemingly different, but is of the same nature.

Needless to say, this bijection provides a correspondence between
down-up permutations and even increasing trees, without using the language of the up-down labelings.

\vskip 6mm \noindent{\large\bf Acknowledgment.} We wish to thank the referee for helpful comments and suggestions.
 This work was supported
by the National
Science Foundation of China.

\end{document}